\documentclass[11pt]{article}
\usepackage[margin=1in]{geometry}
\usepackage{amsmath,amsfonts,amsthm,amssymb}
\usepackage{subcaption}
\usepackage{graphicx}
\graphicspath{{figures/}}
\usepackage{color}
\usepackage{url}
\usepackage{soul}

\usepackage{verbatim}
\usepackage{hyperref}
\DeclareMathOperator{\Res}{Res}
\DeclareMathOperator{\diag}{diag}

\title{Regularity of Solutions for the Peridynamics Equation on Periodic Distributions}
\author{Thinh Dang, Bacim Alali, and Nathan Albin\\
\
\footnotesize{Department of Mathematics, Kansas State University, Manhattan, KS }}
\date{}

\newcommand{\Lper}{\mathcal{L}}
\newcommand{\uper}{\mathbf{u}}
\newcommand{\T}{\mathbb{T}}

\newtheorem{definition}{Definition}
\newtheorem{theorem}{Theorem}
\newtheorem{lemma}{Lemma\textbf{}}
\newtheorem{proposition}{Proposition\textbf{}}
\newtheorem{cor}{Corollary}

\theoremstyle{remark}
\newtheorem{remark}{Remark}
\begin{document}
\title{Regularity of Solutions for  Peridynamics Equilibrium and Evolution Equations on Periodic Distributions}
          \author{Thinh Dang	\thanks{
          tdang@ksu.edu, Department of Mathematics, Kansas State University, Manhattan, KS.}
          \and Bacim Alali \thanks{bacimalali@math.ksu.edu, Department of Mathematics, Kansas State University, Manhattan, KS.}
          \and  Nathan Albin \thanks{albin@ksu.edu, Department of Mathematics, Kansas State University, Manhattan, KS.}
          \thanks{This project is based upon work supported by the National Science Foundation under Grant No. 2108588.}}          
\maketitle

\begin{abstract}
    Results on the peridynamics equilibrium and evolution equations over the space of periodic vector-distributions in multi-spatial dimensions are presented. The associated operator considered is  the linear state-based peridynamic operator for a homogeneous material. Results for weakly singular (integrable) as well as singular  integral kernels are developed. 
    The asymptotic behavior of the eigenvalues of the peridynamic operator's Fourier multipliers and eigenvalues are characterized explicitly in terms of the nonlocality (peridynamic horizon), the integral kernel singularity, and the spatial dimension. 
    We build on the asymptotic analysis to develop   regularity of solutions results for the peridynamic equilibrium as well as the peridynamic evolution equations over periodic distribution. The regularity results are presented explicitly in terms of the data, the integral kernel singularity, and the spatial dimension. 
    Nonlocal-to-local convergence results are presented for the eigenvalues of the peridynamic operator and for the solutions of the equilibrium and evolution equations. The local limiting behavior is shown for two types of limits as the peridynamic horizon (nonlocality) vanishes or as the integral kernel becomes hyper-singular.
\end{abstract}
\section{Introduction}
In this work, we study the regularity of solutions to the peridynamic  equation
\begin{align*}
    \begin{cases}
    \mathbf{u}_{tt}(x,t)&=\Lper \mathbf{u}(x,t) +\mathbf{b}(x),\quad x\in\T^n, \; t>0\\
    \mathbf{u}(x,0)&=\mathbf{f}(x),\quad x\in\T^n\\
    \mathbf{u}_t(x,0)&=\mathbf{g}(x),\quad x\in\T^n,
    \end{cases}
\end{align*}
and the peridynamic equilibrium equation
\[
\Lper \mathbf{u}(x) =\mathbf{b}(x),\quad x\in\T^n.
\]
Here $\T^n$ is the periodic torus in $\mathbb{R}^n$, with $n\ge 1$,  $\mathbf{u}$ is the displacement field, $\mathbf{f}:\T^n\to\mathbb{R}^n$ and $\mathbf{g}:\T^n\to\mathbb{R}^n$ are the initial displacement and the initial velocity, respectively, and $\mathbf{b}:\T^n\to\mathbb{R}^n$ is the external force field. The operator $\Lper$ is the state-based linear peridynamic operator defined by \eqref{eq:explicit_operator}, see Section~ \ref{sec:Peri} for more details about this operator.

There have been some results on the regularity of solutions for peridynamics and nonlocal integral equations of the peridynamic type. Existence and uniqueness for the peridynamic equilibrium equation on bounded domains was proved in \cite{mengesha2014nonlocal}. Existence and uniqueness for the  peridynamic evolution equation was proved in \cite{emmrichperidynamic} and in \cite{alali2012multiscale}.  The works in \cite{foss2018existence, foss2016differentiability} showed some regularity results for the nonlocal Poisson  equation on bounded domains for  integrable kernels. 

The works in \cite{alimov2014problems,alimov2019solvability,alimov2020solvability,alimov2023hypersingular} provide regularity results for the bond-based peridynamics evolution equation for different cases related to the integral kernels and the spatial dimension.
 In this work, we improve on previous results and provide a comprehensive description of the regularity of solutions for  peridynamic equations   over the space of periodic distributions in any spatial dimension. The results are presented for  linear state-based peridynamics equilibrium and evolution equations. The analysis results apply to bond-based peridynamics as a special case of the more general state-peridynamics.
 
The approach followed in this work builds on previous works on scalar peridynamic/nonlocal equations given in \cite{alali2021fourier,ilyas,dang2024}. 
In  \cite{alali2021fourier}, the Fourier multipliers of the nonlocal Laplacian were derived and represented in an integral form as well as in terms of hypergeometric functions, which was key to characterize the multipliers asymptotic behavior and then to develop regularity of solutions analysis for the periodic nonlocal Poisson  equation. 

This work builds also on the work in \cite{alali2022}, which studied the tensor multipliers of the state-based peridynamic operator. It was shown that the multipliers and the eigenvalues of the peridynamic operator have integral representations as well as  hypergeometric representations. These representations are central to the current work and allow for   studying the asymptotic behavior of the tensor multipliers and the eigenvalues and their nonlocal-to-local limiting behavior.

In this work, we generalize the approach in \cite{alali2021fourier}, building on the explicit representations of the multipliers in \cite{alali2022}, to characterize the asymptotic behavior of the state-based peridynamic operator. We then utilize the asymptotic analysis to establish  spatial and temporal regularity of solutions for the equilibrium peridynamic equation as well as for the peridynamic evolution equation. Furthermore, we  study the convergence of the solutions of peridynamic equations to their local counterparts governed by the Navier operator. The nonlocal-to-local convergence of solutions analysis is shown for two types of limiting behavior as the peridynamic horizon (nonlocality) vanishes or as the integral kernel becomes hyper-singular.

The organization  and the main contributions of this article are described  as follows. For completeness of the presentation,  we  provide a review in  Section \ref{sec:Peri} of some results on the peridynamic operator, which appeared in \cite{alali2022}.   Then, we present results on the asymptotic behavior of the peridynamic operator's Fourier multipliers and eigenvalues in Section \ref{sec:eigen-assymp}. In Section \ref{sec:Operator-convergence}, we present nonlocal-to-local convergence results for the peridynamic operator as well as for its eigenvalues. Section \ref{sec:Sobolev} provides an overview of the Sobolev space  of periodic distributions, which will be the set-up space for the equations  considered in subsequent sections. 
The peridynamic equilibrium equation's regularity of solutions and their nonlocal-to-local convergence are developed in Section \ref{sec:Equilibrium}. 
The peridynamic evolution equation is studied in Section \ref{sec:temp-homo} and Section \ref{sec:temp-forcing}.
We study the spatial regularity of these equations over the space of periodic vector-distributions. 
The temporal regularity of solutions, presented in Section \ref{sec:temp-homo} and Section \ref{sec:temp-forcing}, is developed using Gateaux derivative, when the solutions are distributions, and using classical derivative, when the solutions are regular $(L^2(\T^n))^n$ functions.
In Sections \ref{sec:Equilibrium}, \ref{sec:Conv-homo}, and \ref{sec:Conv-forcing}, the convergence of solutions to their local counterparts are shown under two types of limits as the spatial nonlocality vanishes or as the integral kernel becomes hyper-singular.


\section{Linear peridynamic operators}\label{sec:Peri}
 For an homogeneous isotropic solid, the linear peridynamic operator has the form
\begin{align*}
    \Lper\uper(x)&=\rho\int_{\Omega}\frac{\gamma(\|y-x\|)}{\|y-x\|^2}(y-x)\otimes(y-x)(\uper(y)-\uper(x))dy\\
    &+\rho\prime\int_{\Omega}\int_{\Omega}\gamma(\|y-x\|)\gamma(\|z-x\|)(y-x)\otimes(z-x)(\uper(z)-\uper(x))dzdy\\
    &+\rho\prime\int_{\Omega}\int_{\Omega}\gamma(\|y-x\|)\gamma(\|z-y\|)(y-x)\otimes(z-y)(\uper(z)-\uper(y))dzdy,
\end{align*}
where $\gamma$ is a scalar field, and $\rho,\rho\prime$
are scaling constants including the material properties. In the case $\Omega=\mathbb{R}^n$, due to the symmetry, the formula reduces to
\begin{align*}
    \Lper\uper(x)&=\rho\int_{\mathbb{R}^n}\frac{\gamma(\|y-x\|)}{\|y-x\|^2}(y-x)\otimes(y-x)(\uper(y)-\uper(x))dy\\
    &+\rho\prime\int_{\mathbb{R}^n}\int_{\mathbb{R}^n}\gamma(\|y-x\|)\gamma(\|z-y\|)(y-x)\otimes(z-y)\uper(z)dzdy.
\end{align*}
In this work, we will focus  on radially symmetric kernels with compact support of the form
\begin{align*}
    \gamma(\|y-x\|)=c^{\delta,\beta}\frac{1}{\|y-x\|^\beta}\chi_{B_\delta(x)}(y),
\end{align*}
where $c^{\delta,\beta}$ is given by
\begin{align*}
    c^{\delta,\beta}:=\frac{2(n+2-\beta)\Gamma\left(\frac{n}{2}+1\right)}{\pi^{n/2}\delta^{n+2-\beta}}.
\end{align*}
Here, $\chi_{B_\delta(x)}$ denotes the indicator function of the ball centered at $x$ with radius $\delta>0$,  and the kernel exponent $\beta$ satisfies $\beta<n+2$. The linear peridynamic operator, parameterized by the nonlocality parameter $\delta$ and the integral kernel exponent $\beta$ can then be represented as
\begin{align}
\label{eq:explicit_operator}   
\nonumber\Lper^{\delta,\beta}\uper(x)&=(n+2)\mu c^{\delta,\beta}\int_{B_\delta(x)}\frac{(y-x)\otimes(y-x)}{\|y-x\|^{\beta+2}}(\uper(y)-\uper(x))dy\\
    &+(\lambda^*-\mu)\frac{(c^{\delta,\beta})^2}{4}\int_{B_\delta(x)}\int_{B_\delta(y)}\frac{y-x}{\|y-x\|^\beta}\otimes\frac{z-y}{\|z-y\|^\beta}\uper(z)dzdy,
\end{align}
where $\mu$ and $\lambda^*$ are Lamé parameters. (Here we use $\lambda^*$ to denote the second Lamé parameter and save $\lambda$ to denote the eigenvalues later).

\noindent
For convenience, we will decompose $\Lper^{\delta,\beta}$ into the sum of two operators $\Lper_b$ and $\Lper_s$, where these operators refer to the first and second term in the summation respectively. After changing variables, we can rewrite those operators as following
\begin{align}\label{eq:Lb}
    \Lper_b\uper(x)=(n+2)\mu c^{\delta,\beta}\int_{B_\delta(0)}\frac{w\otimes w}{\|w\|^{\beta+2}}(\uper(x+w)-\uper(x))dw,
\end{align}
and
\begin{align}\label{eq:Ls}
    \Lper_s\uper(x)=(\lambda^*-\mu)\frac{(c^{\delta,\beta})^2}{4}\int_{B_\delta(0)}\int_{B_\delta(0)}\frac{w}{\|w\|^\beta}\otimes\frac{q}{\|q\|^\beta}\uper(x+w+q)dqdw.
\end{align}
\noindent
We will be using the following results , Theorem \ref{thm:multipleirs} - Theorem \ref{thm:eigen_integral}, on the multipliers of $\Lper_b,\Lper_s$, and $\Lper^{\delta,\beta}$, which were presented in~\cite{alali2022}.
\begin{theorem}\label{thm:multipleirs}
Expressing $\uper$ through its Fourier transform as
\begin{align*}
    \uper(x)=\frac{1}{(2\pi)^n}\int_{\mathbb{R}^n}\hat{\uper}(\nu)e^{i\nu\cdot x}d\nu,
\end{align*}
we have the representations
\begin{align*}
    &\Lper_b\uper(x)=\frac{1}{(2\pi)^n}\int_{\mathbb{R}^n}M_b(\nu)\hat{\uper}(\nu)e^{i\nu\cdot x}d\nu,\\
    &\Lper_s\uper(x)=\frac{1}{(2\pi)^n}\int_{\mathbb{R}^n}M_s(\nu)\hat{\uper}(\nu)e^{i\nu\cdot x}d\nu,
\end{align*}
where
\begin{align*}
    M_b(\nu)=(n+2)\mu c^{\delta,\beta}\int_{B_\delta(0)}\frac{w\otimes w}{\|w\|^{\beta+n}}(\cos(\nu\cdot x)-1)dw,
\end{align*}
and
\begin{align*}
    M_s(\nu)=-(\lambda^*-\mu)\frac{(c^{\delta,\beta})^2}{4}\left(\int_{B_\delta(0)}\frac{w}{\|w\|^\beta}\sin(\nu\cdot w)dw\right)\otimes\left(\int_{B_\delta(0)}\frac{w}{\|w\|^\beta}\sin(\nu\cdot w)dw\right).
\end{align*}
In addition, the multipliers $M^{\delta,\beta}$ of $\Lper^{\delta,\beta}$ is given by $M^{\delta,\beta}=M_b+M_s$, which satisfies $\widehat{\Lper^{\delta,\beta}\uper}=M^{\delta,\beta}\hat{\uper}$. 
\end{theorem}

\begin{theorem}\label{thm:hyper_Mb}
The hypergeometric representation of $M_b(\nu)$ is given by
\begin{align*}
    M_b(\nu)=\alpha_{b1}(\nu)I+\alpha_{b2}(\nu)\nu\otimes\nu,
\end{align*}
where $I$ denotes the identity matrix and 
\begin{align*}
    &\alpha_{b1}(\nu)=-\mu\|\nu\|^2{}_2F_3\left(1,\frac{n+2-\beta}{2};2,\frac{n+4}{2},\frac{n+4-\beta}{2};-\frac{1}{4}\|\nu\|^2\delta^2\right),\\
    &\alpha_{b2}(\nu)=-2\mu\;{}_1F_2\left(\frac{n+2-\beta}{2};\frac{n+4}{2},\frac{n+4-\beta}{2};-\frac{1}{4}\|\nu\|^2\delta^2\right).
\end{align*}
\end{theorem}

\begin{theorem}\label{thm:hyper_Ms}
The hypergeometric representation of $M_s(\nu)$ is given by
\begin{align*}
    M_s(\nu)=\alpha_s(\nu)\nu\otimes\nu,
\end{align*}
where
\begin{align*}
    \alpha_s(\nu)=-(\lambda^*-\mu)\;{}_1F_2\left(\frac{n+2-\beta}{2};\frac{n+2}{2},\frac{n+4-\beta}{2};-\frac{1}{4}\|\nu\|^2\delta^2\right)^2.
\end{align*}
\end{theorem}

\begin{theorem}\label{thm:hyper_M}
The hypergeometric representation of the multiplier $M^{\delta,\beta}(\nu)$ is given by
\begin{align*}
    M^{\delta,\beta}(\nu)=\alpha_{b1}(\nu)I+(\alpha_{b2}(\nu)+\alpha_s(\nu))\nu\otimes \nu.
\end{align*}
This shows that $M^{\delta,\beta}(\nu)$ is a real symmetric matrix. Furthermore, $\nu$ is an eigenvector of $M^{\delta,\beta}(\nu)$, that is $$M^{\delta,\beta}(\nu)\nu=\lambda_1(\nu)\nu,$$ where
\begin{align*}
    \lambda_1(\nu)&=\alpha_{b1}(\nu)+(\alpha_{b2}(\nu)+\alpha_s(\nu))\|\nu\|^2\\
    &=-\|\nu\|^2\left(3\mu\;{}_3F_4\left(1,\frac{5}{2},\frac{n+2-\beta}{2};2,\frac{3}{2},\frac{n+4}{2},\frac{n+4-\beta}{2};-\frac{1}{4}\|\nu\|^2\delta^2\right)\right.\\
    &\left. \qquad\qquad +(\lambda^*-\mu)\;{}_1F_2\left(\frac{n+2-\beta}{2};\frac{n+2}{2},\frac{n+4-\beta}{2};-\frac{1}{4}\|\nu\|^2\delta^2\right)^2\right).
\end{align*}
The other $n-1$ eigenvectors are orthogonal to $\nu$. Denote any vector in $\mathbb{R}^n$ orthogonal to $\nu$ as $\nu^\bot$, then
\begin{align*}
    M^{\delta,\beta}(\nu)\nu^\bot=\lambda_2(\nu)\nu^\bot,
\end{align*}
where
\begin{align*}
    \lambda_2(\nu)=\alpha_{b1}(\nu)=-\mu\|\nu\|^2\;{}_2F_3\left(1,\frac{n+2-\beta}{2};2,\frac{n+4}{2},\frac{n+4-\beta}{2};-\frac{1}{4}\|\nu\|^2\delta^2\right).
\end{align*}
\end{theorem}

\begin{theorem}\label{thm:eigen_integral}
For any $\nu$ in $\mathbb{R}^n$, the eigenvalues $\lambda_1(\nu)$ and $\lambda_2(\nu)$ of $M^{\delta,\beta}(\nu)$, corresponding to eigenvectors that are parallel to $\nu$ and orthogonal to $\nu$, respectively, have the following integral representations
\begin{equation}
    \begin{aligned}\label{lambda1}
        \lambda_1(\nu)=(n+2)\mu c^{\delta,\beta}&\int_{B_\delta(0)}\frac{(\nu\cdot w)^2}{\|\nu\|^2\|w\|^{\beta+2}}(\cos(\nu\cdot w)-1)dw\\
    &-(\lambda^*-\mu)\left(\frac{c^{\delta,\beta}}{2}\int_{B_\delta(0)}\frac{\nu\cdot w}{\|\nu\|\|w\|^\beta}\sin(\nu\cdot w)dw\right)^2,
    \end{aligned}
\end{equation}
and
\begin{align}\label{lambda2}
    \lambda_2(\nu)=(n+2)\mu c^{\delta,\beta}\int_{B_\delta(0)}\frac{\nu\cdot w}{\|\nu\|^2 \|w\|^{\beta+2}}(\sin(\nu\cdot w)-\nu\cdot w)dw.
\end{align}
\end{theorem}
\noindent
From \eqref{lambda2}, one can claim that $\lambda_2(\nu)<0$ for any nonzero $\nu$, while from \eqref{lambda1}, the sign of $\lambda_1(\nu)$ depends on the values of parameters $\lambda^*$ and $\mu$. Throughout the remaining of this paper, we will assume that some additional conditions are imposed on $\lambda^*$ and $\mu$ such that $\lambda_1(\nu)$ is also negative for any nonzero $\nu$. With this assumption, we can easily show that the multiplier matrices $M^{\delta,\beta}(\nu)$ are invertible for any nonzero $\nu$. This is obvious due to the fact that \begin{align*}
    \det M^{\delta,\beta}(\nu)=\lambda_1(\nu)\lambda_2^{n-1}(\nu).
\end{align*}

\subsection{Eigenvalues' asymptotic analysis }\label{sec:eigen-assymp}
The asymptotic analysis for the eigenvalues of the peridynamic operator $\Lper$, provided in this section, is key to studying the regularity of solutions of peridynamic equations in subsequent sections.

Following a similar approach to that  used in \cite{alali2021fourier}, we analyze the asymptotic behavior of the eigenvalues $\lambda_1(\nu)$ and $\lambda_2(\nu)$ as $\|\nu\|\to\infty$. The asymptotic behavior of $\lambda_2(\nu)$ is summarized in Theorem \ref{thm:asymptotics_2} below.

\begin{theorem}\label{thm:asymptotics_2}
 Let $n\ge 1$, $\delta>0$ and 
 $\beta <n+2$.  Then, as $\|\nu\|\to\infty$,
\begin{equation*}
\lambda_2(\nu) \sim 
\begin{cases}\displaystyle
-\mu\frac{2(n+2-\beta)(n+2)}{\delta^2(n-\beta)}-\mu\frac{\Gamma(\frac{n+4}{2})\Gamma(\frac{n+4-\beta}{2})}{\frac{\beta-n}{2}\Gamma(\frac{\beta+2}{2})}\left(\frac{2}{\delta}\right)^{n+2-\beta}\|\nu\|^{\beta-n}
&\text{if $\beta\ne n$},\\
-\frac{2\mu}{\delta^2}(n+2)\left(
2\log\|\nu\|+
\log\left(\frac{\delta^2}{4}\right)+\gamma-\psi(\frac{n+2}{2})\right)
&\text{if $\beta = n$},
\end{cases}
\end{equation*}
where $\gamma$ is Euler's constant and $\psi$ is the digamma function.
\end{theorem}

\begin{proof}
Let $a=\frac{n+2-\beta}{2}$, $b=\frac{n+4}{2}$ and $z=\frac{\delta\|\nu\|}{2}$ then
\begin{equation}\label{eq:lambda_2}
    \lambda_2(\nu)=-\mu\|\nu\|^2{}_2F_3(1,a;2,b,a+1;-z^2).
\end{equation}
For large $|z|$, ~\cite[Eq.~(16.11.8)]{NIST:DLMF} states that 
\begin{equation*}
_2F_3(1,a;2,b,a+1;-z^2) 
\sim a\Gamma(b)\left(H_{2,3}(z^2) + E_{2,3}(z^2e^{-i\pi}) + E_{2,3}(z^2e^{i\pi})\right),
\end{equation*}
where $E_{2,3}$ and $H_{2,3}$ are formal series defined in~\cite[Eq.~(16.11.1)]{NIST:DLMF} and~\cite[Eq.~(16.11.2)]{NIST:DLMF} respectively. From those definitions, it yields that
\begin{equation*}
E_{2,3}(z^2e^{\pm i\pi}) = (2\pi)^{-1/2}2^{b+1}e^{\pm 2iz}
\sum_{k=0}^{\infty}c_k(\pm 2iz)^{-\frac{2b+3}{2}-k}.
\end{equation*}
Since $2b=n+4$ and $c_0=1$, the $E_{2,3}$ terms decay asymptotically like $|z|^{-\frac{n+7}{2}}$, and do not contribute the the asymptotic behavior described in the theorem.
\noindent
For the remaining term that relates to the $H_{2,3}$ function, from the remark below~\cite[Eq.~(16.11.5)]{NIST:DLMF}, $H_{2,3}$ can be recognized as the sum of the residues of certain poles of the integrand in~\cite[Eq.~(16.5.1)]{NIST:DLMF}. The integrand for this setting is as following:
\begin{equation*}
f(s) := \frac{\Gamma(1+s)\Gamma(a+s)}{\Gamma(2+s)\Gamma(b+s)\Gamma(a+1+s)}
\Gamma(-s)(z^2)^s = \frac{\Gamma(-s)}{(s+1)(s+a)\Gamma(b+s)}z^{2s}.
\end{equation*}
The two poles of interest are at $s=-1$ and $s=-a$.  The restriction that $\beta\notin\{n+2,n+4,n+6,\ldots\}$ ensures that $-a$ is not a nonnegative integer and, therefore, that $s=-a$ is not a pole of $\Gamma(-s)$.  Thus, there are only two cases to consider, either $\beta\ne n$, which implies that $s=-1$ and $s=-a$ are distinct simple poles of $f$, or $\beta=n$, yielding a double pole at $s=-1$.

\noindent
In the case where $\beta\ne n$, we have
\begin{equation*}
H_{2,3}(z^2) = \Res(f,-1) + \Res(f,-a)
= \frac{1}{(a-1)\Gamma(b-1)}z^{-2} + \frac{\Gamma(a)}{(1-a)\Gamma(b-a)}z^{-2a}.
\end{equation*}
When $\beta=n$, the double pole makes the residue more complicated (and gives rise to the logarithmic term):
\begin{equation*}
H_{2,3}(z^2)=\Res(f,-1) = \left.\frac{d}{ds}\left(\frac{\Gamma(-s)}{\Gamma(b+s)}z^{2s}\right)\right|_{s=-1} = \frac{2\log z-\psi(b-1)-\psi(1)}{\Gamma(b-1)}z^{-2}.
\end{equation*}
Finally, using~\eqref{eq:lambda_2} and substituting completes the theorem.
\end{proof}

For the asymptotic behavior of $\lambda_1(\nu)$, we first decompose $\lambda_1(\nu)$ into 
\begin{equation*}
    \lambda_1(\nu)=\lambda_{1,1}(\nu)+\lambda_{1,2}(\nu),
\end{equation*}
where 
\begin{align*}
    &\lambda_{1,1}(\nu)=-3\mu\|\nu\|^2\: {}_3F_4\left(1,\frac{5}{2},\frac{n+2-\beta}{2};2,\frac{3}{2},\frac{n+4}{2},\frac{n+4-\beta}{2};-\frac{1}{4}\|\nu\|^2\delta^2\right),\\
    &\lambda_{1,2}(\nu)=
    -\|\nu\|^2(\lambda^*-\mu)\;{}_1F_2\left(\frac{n+2-\beta}{2};\frac{n+2}{2},\frac{n+4-\beta}{2};-\frac{1}{4}\|\nu\|^2\delta^2\right)^2,
\end{align*}
and then obtain results for $\lambda_{1,1}(\nu)$ in Theorem \ref{thm:asymptotics_{1,1}}  and for $\lambda_{1,2}(\nu)$ in Theorem \ref{thm:asymptotics_{1,2}}.

\begin{theorem}\label{thm:asymptotics_{1,1}}
Let $n\ge 1$, $\delta>0$ and $\beta<n+2$.  Then, as $\|\nu\|\to\infty$,
\begin{equation*}
\lambda_{1,1}(\nu) \sim 
\begin{cases}\displaystyle
-\mu\frac{2(n+2-\beta)(n+2)}{\delta^2(n-\beta)}-2\mu\frac{n-\beta-1}{n-\beta}\frac{\Gamma(\frac{n+4}{2})\Gamma(\frac{n+4-\beta}{2})}{\Gamma(\frac{\beta+2}{2})}\left(\frac{2}{\delta}\right)^{n+2-\beta}\|\nu\|^{\beta-n}
&\text{if $\beta\ne n$},\\
-\frac{2\mu}{\delta^2}(n+2)\left(
2\log\|\nu\|+
\log\left(\frac{\delta^2}{4}\right)+\gamma+2-\psi(\frac{n+2}{2})\right)
&\text{if $\beta = n$},
\end{cases}
\end{equation*}
where $\gamma$ is Euler's constant and $\psi$ is the digamma function.
\end{theorem}

\begin{proof}
Let $a=\frac{n+2-\beta}{2}$,$b=\frac{n+4}{2}$ and $z=\frac{\delta\|\nu\|}{2}$ then
\begin{equation}\label{eq:lambda_{1,1}}
    \lambda_{1,1}(\nu)=-3\mu\|\nu\|^2{}_3F_4\left(1,\frac{5}{2},a;2,\frac{3}{2},b,a+1;-z^2\right).
\end{equation}
From the setting,~\cite[Eq.~(16.11.8)]{NIST:DLMF} states that for large $|z|$
\begin{equation*}
{}_3F_4\left(1,\frac{5}{2},a;2,\frac{3}{2},b,a+1;-z^2\right)
\sim\frac{2}{3} a\Gamma(b)\left(H_{3,4}(z^2) + E_{3,4}(z^2e^{-i\pi}) + E_{3,4}(z^2e^{i\pi})\right),
\end{equation*}
where $E_{3,4}$ and $H_{3,4}$ are formal series defined in~\cite[Eq.~(16.11.1)]{NIST:DLMF} and~\cite[Eq.~(16.11.2)]{NIST:DLMF} respectively. From those definitions, it yields that
\begin{equation*}
E_{3,4}(z^2e^{\pm i\pi}) = (2\pi)^{-1/2}2^{b}e^{\pm 2iz}
\sum_{k=0}^{\infty}c_k(\pm 2iz)^{-\frac{2b+1}{2}-k}.
\end{equation*}
Since $2b=n+4$ and $c_0=1$, the $E_{3,4}$ terms decay asymptotically like $|z|^{-\frac{n+5}{2}}$, and do not contribute the the asymptotic behavior described in the theorem.
\noindent
For the remaining term that relates to the $H_{3,4}$ function, from the remark below~\cite[Eq.~(16.11.5)]{NIST:DLMF}, $H_{3,4}$ can be recognized as the sum of the residues of certain poles of the integrand in~\cite[Eq.~(16.5.1)]{NIST:DLMF}. The integrand for this setting is as following:
\begin{equation*}
f(s) := \frac{\Gamma(1+s)\Gamma(\frac{5}{2}+s)\Gamma(a+s)}{\Gamma(2+s)\Gamma(\frac{3}{2}+s)\Gamma(b+s)\Gamma(a+1+s)}
\Gamma(-s)(z^2)^s
=\frac{(\frac{3}{2}+s)\Gamma(-s)}{(s+1)(s+a)\Gamma(b+s)}z^{2s}.
\end{equation*}
The two poles of interest are at $s=-1$ and $s=-a$.  The restriction that $\beta\notin\{n+2,n+4,n+6,\ldots\}$ ensures that $-a$ is not a nonnegative integer and, therefore, that $s=-a$ is not a pole of $\Gamma(-s)$.  Thus, there are only two cases to consider, either $\beta\ne n$, which implies that $s=-1$ and $s=-a$ are distinct simple poles of $f$, or $\beta=n$, yielding a double pole at $s=-1$.

\noindent
In the case where $\beta\ne n$, we have
\begin{equation*}
H_{3,4}(z^2) = \Res(f,-1) + \Res(f,-a)
= \frac{\frac{1}{2}}{(a-1)\Gamma(b-1)}z^{-2} + \frac{(\frac{3}{2}-a)\Gamma(a)}{(1-a)\Gamma(b-a)}z^{-2a}.
\end{equation*}
When $\beta=n$, the double pole makes the residue more complicated (and gives rise to the logarithmic term):
\begin{align*}
H_{3,4}(z^2)=\Res(f,-1) &= \left.\frac{d}{ds}\left(\frac{(\frac{3}{2}+s)\Gamma(-s)}{\Gamma(b+s)}z^{2s}\right)\right|_{s=-1}\\
&=\left.\left[\frac{\Gamma(-s)}{\Gamma(b+s)}z^{2s}+\left(\frac{3}{2}+s\right)\frac{d}{ds}\left(\frac{\Gamma(-s)}{\Gamma(b+s)}z^{2s}\right)\right]\right|_{s=-1}\\
&= \frac{2\log z-\psi(b-1)-\psi(1)+2}{2\Gamma(b-1)}z^{-2}.
\end{align*}
Finally, using~\eqref{eq:lambda_{1,1}} and substituting completes the theorem.
\end{proof}

\begin{theorem}\label{thm:asymptotics_{1,2}}
Let $n\ge 1$, $\delta>0$ and $\beta<n+2$.  Then, as $\|\nu\|\to\infty$,
\begin{equation*}\displaystyle
\lambda_{1,2}(\nu) \sim 
-(\lambda ^*-\mu)\left[\frac{\Gamma(\frac{n+2}{2})\Gamma(\frac{n+4-\beta}{2})}{\Gamma(\frac{\beta}{2})}\right]^2\left(\frac{2}{\delta}\right)^{2(n+2-\beta)}\|\nu\|^{2(\beta-(n+1))}.
\end{equation*}
\end{theorem}

\begin{proof}
Let $a=\frac{n+2-\beta}{2}$,$b=\frac{n+2}{2}$ and $z=\frac{\delta\|\nu\|}{2}$ then
\begin{equation}\label{eq:lambda_{1,2}}
    \lambda_2(\nu)=-\|\nu\|^2(\lambda^*-\mu)\:{}_1F_2(a;b,a+1;-z^2)^2.
\end{equation}
For large $|z|$, ~\cite[Eq.~(16.11.8)]{NIST:DLMF} states that 
\begin{equation*}
_1F_2(a;b,a+1;-z^2) 
\sim a\Gamma(b)\left(H_{1,2}(z^2) + E_{1,2}(z^2e^{-i\pi}) + E_{1,2}(z^2e^{i\pi})\right),
\end{equation*}
where $E_{1,2}$ and $H_{1,2}$ are formal series defined in~\cite[Eq.~(16.11.1)]{NIST:DLMF} and~\cite[Eq.~(16.11.2)]{NIST:DLMF} respectively. From those definitions, it yields that
\begin{equation*}
E_{1,2}(z^2e^{\pm i\pi}) = (2\pi)^{-1/2}2^{b+1}e^{\pm 2iz}
\sum_{k=0}^{\infty}c_k(\pm 2iz)^{-\frac{2b+1}{2}-k}.
\end{equation*}
Since $2b=n+2$ and $c_0=1$, the $E_{1,2}$ terms decay asymptotically like $|z|^{-\frac{n+3}{2}}$, and do not contribute the the asymptotic behavior described in the theorem.
\noindent
For the remaining term that relates to the $H_{1,2}$ function, from the remark below~\cite[Eq.~(16.11.5)]{NIST:DLMF}, $H_{1,2}$ can be recognized as the sum of the residues of certain poles of the integrand in~\cite[Eq.~(16.5.1)]{NIST:DLMF}. The integrand for this setting is as following:
\begin{equation*}
f(s) := \frac{\Gamma(a+s)}{\Gamma(b+s)\Gamma(a+1+s)}
\Gamma(-s)(z^2)^s = \frac{\Gamma(-s)}{(s+a)\Gamma(b+s)}z^{2s}.
\end{equation*}
The only pole of interest is at $s=-a$.  The restriction that $\beta\notin\{n+2,n+4,n+6,\ldots\}$ ensures that $-a$ is not a nonnegative integer and, therefore, that $s=-a$ is not a pole of $\Gamma(-s)$.  Thus
\begin{equation*}
H_{1,2}(z^2) = \Res(f,-a)
= \frac{\Gamma(a)}{\Gamma(b-a)}z^{-2a}.
\end{equation*}
Finally, using~\eqref{eq:lambda_{1,2}} and substituting completes the theorem.
\end{proof}

The asymptotic behavior of $\lambda_1(\nu)$ is summarized in Theorem \ref{thm:asymptotics_1} below.
\begin{theorem}\label{thm:asymptotics_1}
Let $n\ge 1$, $\delta>0$ and $\beta<n+2$.  Then, as $\|\nu\|\to\infty$,
\begin{align*}
\lambda_{1}(\nu) &\sim 
-(\lambda ^*-\mu)\left[\frac{\Gamma(\frac{n+2}{2})\Gamma(\frac{n+4-\beta}{2})}{\Gamma(\frac{\beta}{2})}\right]^2\left(\frac{2}{\delta}\right)^{2(n+2-\beta)}\|\nu\|^{2(\beta-(n+1))}\\
&+\begin{cases}\displaystyle
-\mu\frac{(n+2-\beta)(n+2)}{\delta^2(n-\beta)}-2\mu\frac{n-\beta-1}{n-\beta}\frac{\Gamma(\frac{n+4}{2})\Gamma(\frac{n+4-\beta}{2})}{\Gamma(\frac{\beta+2}{2})}\left(\frac{2}{\delta}\right)^{n+2-\beta}\|\nu\|^{\beta-n}
&\text{if $\beta\ne n$},\\
-\frac{2\mu}{\delta^2}(n+2)\left(
2\log\|\nu\|+
\log\left(\frac{\delta^2}{4}\right)+\gamma+2-\psi(\frac{n+2}{2})\right)
&\text{if $\beta = n$},
\end{cases}
\end{align*}
where $\gamma$ is Euler's constant and $\psi$ is the digamma function.
\end{theorem}

\subsection{Local limits of the multipliers and the peridynamic operator}\label{sec:Operator-convergence}
In this section, we compute the limits of the Fourier multipliers and the peridynamic operator $\Lper$ in \eqref{eq:explicit_operator}. Two types of limits are considered as $\delta\to 0$ with fixed $\beta<n+2$ and as $\beta\to n+2^-$ with fixed $\delta>0$. The results of this section will be used to study the limiting behavior of the solutions of peridynamic equations in subsequent sections. 

Denoting by $\mathcal{N}$ the Navier operator of linear elasticity. For any homogeneous medium, it has the form
\begin{align*}
    \mathcal{N}\mathbf{u}=(\lambda^*+\mu)\nabla(\nabla\cdot\mathbf{u} )+\mu\Delta\mathbf{u}.
\end{align*}
Similar to the peridynamic operator, its multiplier is defined through Fourier transform as $\widehat{\mathcal{N}\mathbf{u}}=M^{\mathcal{N}}\widehat{\mathbf{u}}$ and can be explicitly expressed as
\begin{align*}
    M^{\mathcal{N}}(\nu)=-(\lambda^*+\mu)\nu\otimes\nu-\mu\|\nu\|^2 I.
\end{align*}
The multipliers of the peridynamic operator $M^{\delta,\beta}$ converges to the multiplier of the Navier operator $M^{\mathcal{N}}$ in the following sense (check \cite{alali2022}):
\begin{proposition}\label{prop1}
Let $n\ge 1$, $\delta>0$ and $\beta<n+2$. Then for fixed $\beta$, 
\begin{align*}
    \lim_{\delta\to 0^+}M^{\delta,\beta}(\nu)=M^{\mathcal{N}}(\nu).
\end{align*}
Moreover, for fixed $\delta$,
\begin{align*}
    \lim_{\beta\to n+2^-}M^{\delta,\beta}(\nu)=M^{\mathcal{N}}(\nu).
\end{align*}
\end{proposition}

One corollary of this result is the convergence of the eigenvalues $\lambda_1(\nu)$ and $\lambda_2(\nu)$ of $M^{\delta,\beta}$ to the eigenvalues $\lambda_1^{\mathcal{N}}(\nu)$ and $\lambda_2^{\mathcal{N}}(\nu)$ of $M^{\mathcal{N}}$, respectively.
\begin{cor}
Let $n\ge 1$, $\delta>0$ and $\beta<n+2$. Then for fixed $\beta$,
\begin{align*}
    \lim_{\delta\to 0^+}\lambda_1(\nu)&=-(\lambda^*+\mu)\|\nu\|^2=\lambda_1^{\mathcal{N}}(\nu),\\
    \lim_{\delta\to 0^+}\lambda_2(\nu)&=-\mu\|\nu\|^2=\lambda_2^{\mathcal{N}}(\nu),
\end{align*}
and for fixed $\delta$,
\begin{align*}
    \lim_{\beta\to n+2^-}\lambda_1(\nu)&=-(\lambda^*+\mu)\|\nu\|^2=\lambda_1^{\mathcal{N}}(\nu),\\
    \lim_{\beta\to n+2^-}\lambda_2(\nu)&=-\mu\|\nu\|^2=\lambda_2^{\mathcal{N}}(\nu).
\end{align*}
\end{cor}

The corresponding convergence of the peridynamic operator $\Lper$ to the Navier operator $\mathcal{N}$ as $\delta\to 0$ or as $\beta\to n+2^-$ is given in the following results.

\begin{theorem}
Let $n\ge 1$, $\delta>0$ and $\beta<n+2$. Assume that $\uper(x)\in (C^3(\Omega))^n$, then for fixed $\beta$,
\begin{align*}
    \lim_{\delta\to 0^+}\Lper_b\uper(x)=\mu \Delta \uper (x)+2\mu\nabla(\nabla\cdot \uper )(x),
\end{align*}
and for fixed $\delta$,
\begin{align*}
    \lim_{\beta\to n+2^-}\Lper_b\uper(x)=\mu \Delta \uper (x)+2\mu\nabla(\nabla\cdot \uper )(x).
\end{align*}
\end{theorem}
\begin{proof}
Using Taylor's theorem, $\uper$ is expanded about $x$ as
\begin{align*}
    \uper(x+w)=\uper(x)+\nabla \uper (x):w+\frac{1}{2}\nabla \nabla \uper (x):w\otimes w+R_\sigma(\uper;x,w),
\end{align*}
where the remainder $R_\sigma(\uper;x,w)$ is of the form
\begin{align*}
    R_\sigma(\uper;x,w)=\frac{1}{6}\nabla \nabla \nabla \uper (x+\sigma w):w\otimes w\otimes w,
\end{align*}
for some scalar $\sigma\in [0,1]$. By substituting, \eqref{eq:Lb} becomes
\begin{align*}
    \Lper_b\uper(x)&=(n+2)\mu c^{\delta,\beta} \int_{B_\delta(0)}\frac{w\otimes w\otimes w}{\|w\|^{\beta+2}}dw:\nabla\uper(x)\\
    &+(n+2)\mu c^{\delta,\beta}\int_{B_\delta(0)}\frac{w\otimes w\otimes w\otimes w}{\|w\|^{\beta+2}}dw:\frac{1}{2}\nabla\nabla \uper(x)\\
    &+(n+2)\mu c^{\delta,\beta}\int_{B_\delta(0)}\frac{w\otimes w}{\|w\|^{\beta+2}}R_\sigma(\uper;x,z)dw.
\end{align*}
Due to symmetry, the integral in the first term is identically zero, which makes the whole term vanishes. For the second term, using spherical coordinates, it is straightforward to show that
\begin{equation}\label{Eq1}
    (n+2)c^{\delta,\beta}\int_{B_{\delta}(0)}\frac{w_iw_jw_kw_l}{\|w\|^{\beta+2}}dw=
    \begin{cases}
    6, \quad\text{for}\; i=j=k=l,\\
    2,\quad 
        \text{for}\;i=j\ne k=l,\text{or}\;i=k\ne j=l,\text{or}\;i=l\ne j=k,\\
    0,\quad\text{otherwise}.
    \end{cases}
\end{equation}

Utilizing \eqref{Eq1}, observe that for any component $i$, we have
\begin{align*}
    (n+2) c^{\delta,\beta}\int_{B_\delta(0)}\frac{w_iw_jw_kw_l}{\|w\|^{\beta+2}}dw\frac{1}{2}\frac{\partial^2 \uper_j }{\partial x_l\partial x_k }&=\frac{1}{2}\left[6\frac{\partial^2 \uper_i }{\partial x_i^2}+2\sum_{k\ne i}\frac{\partial^2 \uper_i }{\partial x_k^2}+4\sum_{j\ne i}\frac{\partial^2\uper_j}{\partial x_i\partial x_j}\right]\\
    &=\sum_{k}\frac{\partial^2 \uper_i }{\partial x_k^2}
    +2\sum_{j}\frac{\partial^2\uper_j}{\partial x_i\partial x_j}\\
    &=(\Delta\uper)_i(x)+2\nabla (\nabla\cdot\uper)_i(x). 
\end{align*}
This shows that the second term is $\mu \Delta \uper (x)+2\mu\nabla(\nabla\cdot \uper )(x)$.
Finally, for the third term, we have the following bound
\begin{align*}
    \left|(n+2)\mu c^{\delta,\beta}\int_{B_\delta(0)}\frac{w\otimes w}{\|w\|^{\beta+2}}R_\sigma(\uper;x,z)dw\right|\le M_x\delta (n+2-\beta),
\end{align*}
where $M_x$ is a constant dependent on $x$ only. All these arguments in combination complete the proof.
\end{proof}

\begin{theorem}
Let $n\ge 1$, $\delta>0$ and $\beta<n+2$. Assume that $\uper(x)\in (C^3(\Omega))^n$, then for fixed $\beta$,
\begin{align*}
    \lim_{\delta\to 0^+}\Lper_s\uper(x)=(\lambda^*-\mu)\nabla(\nabla\cdot \uper )(x) \quad \text{for }\beta<n+2,
\end{align*}
and for fixed $\delta$,
\begin{align*}
    \lim_{\beta\to n+2^-}\Lper_s\uper(x)=(\lambda^*-\mu)\nabla(\nabla\cdot \uper )(x) \quad \text{for }\delta>0.
\end{align*}
\end{theorem}

\begin{proof}
Firstly, rewite \eqref{eq:Ls} as
\begin{align}\label{eq:Ls_new}
    \Lper_s\uper(x)=(\lambda^*-\mu)\frac{(c^{\delta,\beta})^2}{4}\int_{B_\delta(0)}\frac{w}{\|w\|^\beta}\left[\int_{B_\delta(0)}\frac{q}{\|q\|^\beta}\cdot \uper(x+w+q)dq\right]dw.
\end{align}
Then $\uper$ is expanded about $x+w$ using Taylor's expansion as
\begin{align*}
    \uper(x+w+q)=\uper(x+w)+\nabla\uper(x+w):q+\frac{1}{2}\nabla\nabla \uper(x+w):q\otimes q+R_{s_1}(\uper;x+w,q),
\end{align*}
where the remainder $R_{s_1}(\uper;x+w,q)$ is of the form
\begin{align*}
    R_{s_1}(\uper;x+w,q)=\frac{1}{6}\nabla \nabla\nabla \uper(x+s_1w):q\otimes q\otimes q,
\end{align*}
for some scalar $s_1\in [0,1]$. 
The inner integral in \eqref{eq:Ls_new} now becomes
\begin{align*}
    \int_{B_\delta(0)}\frac{q}{\|q\|^\beta}\cdot &\uper(x+w+q)dq=\int_{B_\delta(0)}\frac{q}{\|q\|^\beta}dq\cdot \uper(x+w)+\int_{B_\delta(0)}\frac{q\otimes q}{\|q\|^\beta}dq:\nabla
    \uper(x+w)\\
    &\quad+\int_{B_\delta(0)}\frac{q\otimes q\otimes q}{\|q\|^\beta}dq:\frac{1}{2}\nabla\nabla \uper(x+w)+\int_{B_\delta(0)}\frac{q}{\|q\|^\beta}\cdot R_{s_1}(\uper;x+w,q)dq
\end{align*}
Note that due to symmetry, the integral in the first and third term in the expression above vanish. In addition, since
\begin{align*}
    \int_{B_\delta(0)}\frac{q\otimes q}{\|q\|^\beta}dq=\frac{2}{c^{\delta,\beta}}I,
\end{align*}
then
\begin{align}
    \Lper_s\uper(x)&=(\lambda^*-\mu)\frac{c^{\delta,\beta}}{2}\int_{B_\delta(0)}\frac{w}{\|w\|^\beta}\nabla\cdot\uper(x+w)dw\label{Eq2.1}\\
    &+(\lambda^*-\mu)\frac{(c^{\delta,\beta})^2}{4}\int_{B_\delta(0)}\frac{w}{\|w\|^\beta}\left[\int_{B_\delta(0)}\frac{q}{\|q\|^\beta}\cdot R_{s_1}(\uper;x+w,q)dq\right]dw\label{Eq2.2}
\end{align}
Using Taylor's theorem, we have the following expansion
\begin{align*}
    \nabla\cdot\uper(x+w)=\nabla\cdot\uper(x)+\nabla(\nabla\cdot\uper)(x):w+R_{s_2}(\uper;x,w),
\end{align*}
where the remainder $R_{s_2}(\uper;x,w)$ is of the form
\begin{align*}
    R_{s_2}(\uper;x,w)=\frac{1}{2}\nabla\nabla(\nabla\cdot\uper)(x+s_2w):w\otimes w,
\end{align*}
for some scalar $s_2\in [0,1]$. The integral in \eqref{Eq2.1} becomes
\begin{align}
    \int_{B_\delta(0)}\frac{w}{\|w\|^\beta}\nabla\cdot \uper(x+w)dw&=\int_{B_\delta(0)}\frac{w}{\|w\|^\beta}dw\nabla\cdot \uper(x)+\int_{B_\delta(0)}\frac{w\otimes w}{\|w\|^\beta}dw:\nabla(\nabla\cdot\uper)(x)\label{Eq3.1}\\
    &+\int_{B_\delta(0)}\frac{w\otimes w}{\|w\|^\beta}R_{s_2}(\uper;x,w)dw\label{Eq3.2}
\end{align}
In the above expression, the first term in \eqref{Eq3.1} is identically $0$ due to symmetry. About the second term, we have
\begin{align*}
    \int_{B_\delta(0)}\frac{w\otimes w}{\|w\|^\beta}dw:\nabla(\nabla\cdot\uper)(x)=\frac{2}{c^{\delta,\beta}}\nabla(\nabla\cdot\uper)(x).
\end{align*}
We then get from \eqref{Eq2.1}
\begin{align}\label{Eq4}
    \frac{c^{\delta,\beta}}{2}\int_{B_\delta(0)}\nabla\cdot\uper(x+w)\frac{w}{\|w\|^\beta}dw=\nabla(\nabla\cdot\uper)(x)+\frac{c^{\delta,\beta}}{2}\int_{B_\delta(0)}\frac{w\otimes w}{\|w\|^\beta}R_{s_2}(\uper;x,w).
\end{align}
The second expression in \eqref{Eq4} vanishes under the limit of $\delta\to 0^+$ or $\beta \to n+2^-$, where we have used  the bounds
\begin{align*}
    \left|\frac{c^{\delta,\beta}}{2}\int_{B_\delta(0)}\frac{w\otimes w}{\|w\|^\beta}R_{s_2}(\uper;x,w)dw\right|\le M_1\delta(2+n-\beta)
\end{align*}
A similar result applies to \eqref{Eq2.2} as well, from which the statement of the theorem follows.
\end{proof}

\section{The Sobolev space of periodic distributions}\label{sec:Sobolev}
Consider $\Lper^{\delta,\beta}$ as the periodic operator on the torus
\begin{align*}
    \mathbf{T}^n=\prod_{i=1}^n[0,2\pi],\qquad \text{with}\quad i=1,2,\ldots,n.
\end{align*}
The eigenvalues and eigenvectors of the operator $\Lper$ are determined through the multipliers developed in the previous section. The following result, Theorem \ref{thm:L_periodic}, has been proved in~\cite{alali2022}.
\begin{theorem}\label{thm:L_periodic}
Let $\gamma$ be a fixed vector in $\mathbb{R}^n$. For any $k\in\mathbb{Z}^n$, define 
\begin{align*}
    \psi_k(x)=e^{ik\cdot x}\gamma
\end{align*}
Then 
\begin{align*}
    \Lper^{\delta,\beta}\psi_k(x)=M^{\delta,\beta}(k)\psi_k(x)=M_k\psi_k(x).
\end{align*}
\end{theorem}


We define the Sobolev space for periodic vector distributions, which resembles the Sobolev space for periodic scalar distributions defined in \cite{alali2021fourier}.
\begin{definition}
For $q\in\mathbb{R}$, define by $\mathcal{H}^q(\mathbf{T}^n)$ the space of all periodic vector distributions $\mathbf{g}$ on $\mathbf{T}^n$ such that
\begin{align*}
    \|\mathbf{g}\|_{\mathcal{H}^q(\T^n)}^2:=\sum_{k\in\mathbb{Z}^n}(1+\|k\|^2)^q\|\hat{\mathbf{g}}_k\|^2<\infty.
\end{align*}
\end{definition}

From this definition, it shows that if $\mathbf{g}\in\mathcal{H}^q(\T^n)$, then for any nonzero $k\in \mathbb{Z}^n$
\begin{align*}
    \|\widehat{\mathbf{g}}_k\|\le \frac{\|\mathbf{g}\|_{\mathcal{H}^q(\T^n)}}{\|k\|^q}.
\end{align*}

In  subsequent sections, we consider maps from the real line to periodic distributions and use the following notion of derivative for these maps.

\begin{definition}
Suppose $\mathbf{U}:\left[0,\infty\right)\to \mathcal{H}^q(\T^n)$ maps $t\ge 0$ to a periodic vector distribution $\mathbf{U}(t)$ in $\mathcal{H}^q(\T^n)$. The derivative of $\mathbf{U}$ at $t$ is defined in the sense of Gateaux derivative
\begin{align*}
    \mathbf{U}^\prime(t)=\frac{d}{dt}\mathbf{U}(t)=\lim_{h\to 0}\frac{\mathbf{U}(t+h)-\mathbf{U}(t)}{h},
\end{align*}
where the limit is with respect to the norm of $\mathcal{H}^q(\T^n)$ space.\\
We will denote by $C^p([0,\infty), \mathcal{H}^q(\T^n))$ the space of all maps from $[0,\infty)$ to $\mathcal{H}^q(\T^n)$ that is $p$ times differentiable. 
\end{definition}

\begin{lemma}\label{lem:derivative}
Given a map $\mathbf{U}:\left[0,\infty\right)\to \mathcal{H}^q(\T^n)$ and suppose that the derivative $\mathbf{U}^\prime(t)$ exists. Then for any $k\in\mathbb{Z}^n$, the Fourier coefficient $\hat{\mathbf{U}}^\prime_k(t)$ of $\mathbf{U}^\prime(t)$ is the derivative in the classical sense of the Fourier coefficient $\hat{\mathbf{U}}_k(t)$ of $\mathbf{U}(t)$. 
\end{lemma}
\begin{proof}
Represent $\mathbf{U}(t)$ and $\mathbf{U}^\prime(t)$ through their Fourier series
\begin{align*}
    \mathbf{U}(t)=\sum_{k\in\mathbb{Z}^n}\hat{\mathbf{U}}_k(t)e^{i k\cdot x},
\end{align*}
and
\begin{align*}
    \mathbf{U}^\prime(t)=\sum_{k\in\mathbb{Z}^n}\hat{\mathbf{U}}^\prime_k(t)e^{i k\cdot x}.
\end{align*}
From the derivative definition, observe that
\begin{align*}
    0=\lim_{h\to 0}\left\|\frac{\mathbf{U}(t+h)-\mathbf{U}(t)}{h}-\mathbf{U}^\prime(t)\right\|^2_{\mathcal{H}^q(\T^n)}
    =\lim_{h\to 0}\sum_{k\in\mathbb{Z}^n}(1+\|k\|^2)^q\left|\frac{\hat{\mathbf{U}}_k(t+h)-\hat{\mathbf{U}}_k(t)}{h}-\hat{\mathbf{U}}^\prime_k(t)\right|^2.
\end{align*}
The statement implies that
\begin{align*}
    \lim_{h\to 0}\frac{\hat{\mathbf{U}}_k(t+h)-\hat{\mathbf{U}}_k(t)}{h}=\hat{\mathbf{U}}^\prime_k(t).
\end{align*}
for any $k\in \mathbb{Z}^n$. This is the definition for derivative in the classical sense.
\end{proof}

\begin{definition}
Suppose $\mathbf{U}:\left[0,\infty\right)\to \mathcal{H}^q(\T^n)$ such that for any $t\ge 0$, $\mathbf{U}(t)$  has the form
\begin{align*}
    \mathbf{U}(t)=\sum_{k\in\mathbb{Z}^n}\hat{\mathbf{U}}_k(t)e^{ik\cdot x}.
\end{align*}
Let $\Lper^{\delta,\beta}$ and $\mathcal{N}$ be the peridynamic operator and the Navier operator respectively, then define the map $\Lper^{\delta,\beta}\mathbf{U}, \mathcal{N}\mathbf{U}:\left[0,\infty\right)\to \mathcal{H}^q(\T^n)$ through the representation
\begin{align*}
\Lper^{\delta,\beta}\mathbf{U}(t)=\sum_{k\in\mathbb{Z}^n}M_k\hat{\mathbf{U}}_k(t)e^{i k\cdot x},\qquad\text{and}\qquad \mathcal{N}\mathbf{U}(t)=\sum_{k\in\mathbb{Z}^n}M^{\mathcal{N}}_k\hat{\mathbf{U}}_k(t)e^{i k\cdot x},
\end{align*}
for any $t\ge 0$.
\end{definition}

In particular, when $q<0$, the above notion extends the definition of the peridynamic operator $\Lper$ to periodic distributions.

\section{The peridynamic equilibrium equation}\label{sec:Equilibrium}
In this section, we present regularity results for the solutions of the peridynamic equilibrium equation and study their nonlocal-to-local convergence. 

For simplicity, from now on, for any $k\in \mathbb{Z}^n$, we will write $M_k$ and $M_k^{\mathcal{N}}$ in replacement of $M^{\delta,\beta}(k)$ and $M^{\mathcal{N}}(k)$, respectively. Also, since the matrix space is finite dimensional, any matrix norms are equivalent. Throughout this paper, we would use the matrix operator norm that is induced by the vector 2-norm. 

\begin{theorem}
Let $n\ge 1$, $\delta>0$ and $\beta< n+2$. If the vector field $\mathbf{b}\in \mathcal{H}^s(\mathbf{T}^n)$ is such that $\widehat{\mathbf{b}}_0=\mathbf{0}$, then there is uniquely one vector field $\mathbf{u}\in \mathcal{H}^{s^\prime}(\mathbf{T}^n)$ satisfying $\Lper^{\delta,\beta}\mathbf{u}=\mathbf{b}$ and $\widehat{\mathbf{u}}_0=\mathbf{0}$, where $s^\prime=s+\max\{0,\beta-n\}$.
\end{theorem}
\begin{proof}
Represent both $\mathbf{b}$ and $\mathbf{u}$ through their Fourier series respectively and substitute that to the equation $\Lper^{\delta,\beta}\mathbf{u}=\mathbf{b}$ yields
\begin{align*}
    \sum_{k\in \mathbb{Z}^n}\widehat{\mathbf{b}}_k e^{ik\cdot x}=\Lper^{\delta,\beta}\left(\sum_{k\in \mathbb{Z}^n}\widehat{\mathbf{u}}_k e^{ik\cdot x}\right)=\sum_{k\in\mathbf{Z}^n}M_k\widehat{\mathbf{u}}_k e^{ik\cdot x}.
\end{align*}
This shows that for any nonzero $k$, the coefficient vector $\widehat{\mathbf{u}}_k$ is uniquely determined through the formula
\begin{align}\label{F-cof}
    \widehat{\mathbf{u}}_k=(M_k)^{-1}\widehat{\mathbf{b}}_k.
\end{align}
It remains to show that $\mathbf{u}$ is in $\mathcal{H}^{s^\prime}(\mathbf{T}^n)$ where $s^\prime=\max\{0,\beta-n\}$. Observe that
\begin{align*}
    \sum_{k\in\mathbb{Z}^n}(1+\|k\|^2)^{s^\prime}\|\widehat{\mathbf{u}}_k\|^2
    \le\sum_{0\ne k\in\mathbb{Z}^n}(1+\|k\|^2)^{s}\|\widehat{\mathbf{b}}_k\|^2(1+\|k\|^2)^{s^\prime-s}\|(M_k)^{-1}\|^2.
\end{align*}
Since $\mathbf{b}\in \mathcal{H}^s(\mathbf{T}^n)$, it suffices to prove the uniform boundedness of 
\begin{align*}
    (1+\|k\|^2)^{s^\prime-s}\|(M_k)^{-1}\|^2,
\end{align*}
for large value of $\|k\|$.
Since $\lambda_1(k)$ and $\lambda_2(k)$  are eigenvalues of $M_k$, the eigenvalues of $(M_k)^{-1}$ are determined to be $1/\lambda_1(k),1/\lambda_2(k)$. This implies that
\begin{align*}
    \|(M_k)^{-1}\|=\max\left\{\frac{1}{|\lambda_1(k)|},\frac{1}{|\lambda_2(k)|}\right\}=\frac{1}{\min\{|\lambda_1(k)|,|\lambda_2(k)|\}}.
\end{align*}
We then consider two cases. First, when $\beta\le n$, then $s^\prime=s$. Using Theorem~\ref{thm:asymptotics_2} and Theorem~\ref{thm:asymptotics_1}, there exists $r_1>0$ and $C_1>0$ such that $\|(M_k)^{-1}\|\le C_1$ whenever $\|k\|\ge r_1$. This shows that
\begin{align*}
    (1+\|k\|^2)^{s^\prime-s}\|(M_k)^{-1}\|^2\le \frac{1}{C_1^2}.
\end{align*}
For the case when $\beta>n$, we have $s^\prime=s+\beta-n$. Again from Theorem~\ref{thm:asymptotics_2} and Theorem~\ref{thm:asymptotics_1}, there exists $r_2>0$ and $C_2>0$ such that $\|(M_k)^{-1}\|\le \frac{C_2}{\|k\|^{\beta-n}}$ whenever $\|k\|\ge r_2$. This shows that
\begin{align*}
    (1+\|k\|^2)^{s^\prime-s}\|(M_k)^{-1}\|^2\le C_2^2\left(\frac{1+\|k\|^2}{\|k\|^2}\right)^{\beta-n},
\end{align*}
which is uniformly bounded. The proof is now complete.
\end{proof}

For the next two theorems, we provide some results on the convergence of the equilibrium peridynamic equation to the corresponding equilibrium equation of linear elasticity under two types of limits as $\delta\to 0$ with fixed $\beta<n+2$ or as $\beta\to (n+2)^-$ with fixed $\delta> 0$.

\begin{theorem}\label{thm:Equi1}
Let $n\ge 1$, $\delta>0$ and $\beta<n+2$. Suppose $\mathbf{b}\in \mathcal{H}^s(\T^n)$ satisfying $\widehat{\mathbf{b}}_0=\mathbf{0}$. Let $\mathbf{u}^{\mathcal{N}}\in \mathcal{H}^{s+2}(\T^n)$ be the solution of the Navier equation $\mathcal{N}\mathbf{u}^{\mathcal{N}}=\mathbf{b}$ with $\widehat{\mathbf{u}}_0^{\mathcal{N}}=\mathbf{0}$ and let $\mathbf{u}^{\delta,\beta}\in \mathcal{H}^{s'}(\T^n)$ be the solution to the equation $\mathcal{L}^{\delta,\beta}\mathbf{u}=\mathbf{b}$ with $\widehat{\mathbf{u}}_0^{\delta,\beta}=\mathbf{0}$ where $s'=s+\max\{0,\beta-n\}$. Then
\begin{align*}
    \lim_{\delta\to 0^+}\mathbf{u}^{\delta,\beta}=\mathbf{u}^{\mathcal{N}}, \quad \text{in } \mathcal{H}^{s'}(\T^n).
\end{align*}
\end{theorem}

\begin{proof}
From \ref{F-cof}, we see that
\begin{align*}
    \|\mathbf{u}^{\delta,\beta}-\mathbf{u}^{\mathcal{N}}\|^2_{\mathcal{H}^{s'}(\T^n)}&=\sum_{0\ne k \in \mathbb{Z}^n}(1+\|k\|^2)^{s'}\|\widehat{\mathbf{u}}_k^{\delta,\beta}-\widehat{\mathbf{u}}_k^{\mathcal{N}}\|^2\\
    &\le\sum_{0\ne k \in \mathbb{Z}^n}(1+\|k\|^2)^{s'}\|(M_k)^{-1}-(M_k^{\mathcal{N}})^{-1}\|^2\|\widehat{\mathbf{b}}_k\|^2.
\end{align*}
Each term in the summation converges to $0$ under the limit using  Proposition \ref{prop1}. The proof is then complete if we can pass the limit inside the summation. Since $\mathbf{b}\in \mathcal{H}^s(\T^n)$, this is feasible provided that $(1+\|k\|^2)^{s'-s}\|(M_k)^{-1}-(M_k^{\mathcal{N}})^{-1}\|^2$ or equivalently, $\|k\|^{s'-s}\|(M_k)^{-1}-(M_k^{\mathcal{N}})^{-1}\|$ is uniformly bounded for large value of $\|k\|$. We note that 
\begin{align*}
    \|(M_k)^{-1}-(M_k^{\mathcal{N}})^{-1}\|=\frac{1}{\min\{|\lambda_1(k)-\lambda_1^{\mathcal{N}}(k)|,|\lambda_2(k)-\lambda_2^{\mathcal{N}}(k)|\}}.
\end{align*}
Using Theorem \ref{thm:asymptotics_2} and Theorem \ref{thm:asymptotics_1}, one can show that there exists $r, C_1, C_2, C_3>0$ such that whenever $\|k\|\ge r$
\begin{align*}
    |\lambda_1(k)|,|\lambda_2(k)|\le\begin{cases}
    C_1&\text{if }\beta<n\\
    C_2\log\|k\|&\text{if }\beta=n\\
    C_3\|k\|^{\beta-n}&\text{if }n<\beta<n+2
    \end{cases}
\end{align*}
Since $|\lambda_1^{\mathcal{N}}(k)|=|\lambda^*+\mu|\|k\|^2$ and $|\lambda_2^{\mathcal{N}}(k)|=\mu\|k\|^2$, respectively, there is $C>0$ such that for large value of $\|k\|$
\begin{align*}
    \|(M_k)^{-1}-(M_k^{\mathcal{N}})^{-1}\|\le \frac{C}{\|k\|^2}
\end{align*}
This along with the fact that
\begin{align*}
    \|k\|^{s'-s}=\|k\|^{\max\{0,\beta-n\}}\le \|k\|^2
\end{align*}
yields the desired uniform bound.
\end{proof}

\begin{theorem}
Let $n\ge 1$, $\delta>0$ and $\beta<n+2$. Suppose $\mathbf{b}\in \mathcal{H}^s(\T^n)$ satisfying $\widehat{\mathbf{b}}_0=\mathbf{0}$. Let $\mathbf{u}^{\mathcal{N}}\in \mathcal{H}^{s+2}(\T^n)$ be the solution the the equation $\mathcal{N}\mathbf{u}=\mathbf{b}$ with $\widehat{\mathbf{u}}_0^{\mathcal{N}}=\mathbf{0}$ and let $\mathbf{u}^{\delta,\beta}\in H^{s'}(\T^n)$ be the solution to the equation $\mathcal{L}^{\delta,\beta}\mathbf{u}=\mathbf{b}$ with $\widehat{\mathbf{u}}_0^{\delta,\beta}=\mathbf{0}$ where $s'=s+\max\{0,\beta-n\}$ . Then, for any $0<\epsilon<2$, 
\begin{align*}
    \lim_{\beta\to n+2^-}\mathbf{u}^{\delta,\beta}=\mathbf{u}^{\mathcal{N}}, \quad \text{in } \mathcal{H}^{s_0}(\T^n),
\end{align*}
where $s_0=s+2-\epsilon$.
\end{theorem}
\begin{proof}
Fix an $\epsilon\in(0,2)$. Then for any $n+2-\epsilon<\beta<n+2$, both $s+2$ and $s'$ are larger than $s_0$, hence the limit makes sense. Similar to Theorem \ref{thm:Equi1}, we have that
\begin{align*}
    \|\mathbf{u}^{\delta,\beta}-\mathbf{u}^{\mathcal{N}}\|^2_{\mathcal{H}^{s_0}(\T^n)}&=\sum_{0\ne k \in \mathbb{Z}^n}(1+\|k\|^2)^{s_0}\|\widehat{\mathbf{u}}_k^{\delta,\beta}-\widehat{\mathbf{u}}_k^{\mathcal{N}}\|^2\\
    &\le\sum_{0\ne k \in \mathbb{Z}^n}(1+\|k\|^2)^{s_0}\|(M_k)^{-1}-(M_k^{\mathcal{N}})^{-1}\|^2\|\widehat{\mathbf{b}}_k\|^2.
\end{align*}
The proof is established provided the uniform boundedness of 
\begin{align*}
    \|k\|^{2-\epsilon}\|(M_k)^{-1}-(M_k^{\mathcal{N}})^{-1}\|
\end{align*}
Again, there is $C>0$ such that
\begin{align*}
    \|(M_k)^{-1}-(M_k^{\mathcal{N}})^{-1}\|\le \frac{C}{\|k\|^2}
\end{align*}
for large $\|k\|$. This yields that
\begin{align*}
    \|k\|^{2-\epsilon}\|(M_k)^{-1}-(M_k^{\mathcal{N}})^{-1}\|\le \frac{C}{\|k\|^\epsilon},
\end{align*}
which is uniformly bounded.
\end{proof}

\section{The  peridynamic initial-value problem}\label{sec:Peri-homo}
In this section, we focus on the following peridynamic equation with initial data and zero forcing term
\begin{align}\label{eq:Homo-peri}
    \begin{cases}
    \mathbf{u}_{tt}(x,t)&=\Lper \mathbf{u}(x,t),\quad x\in\T^n, \; t>0,\\
    \mathbf{u}(x,0)&=\mathbf{f}(x),\\
    \mathbf{u}_t(x,0)&=\mathbf{g}(x).
    \end{cases}
\end{align}
In order to study the existence, uniqueness, and regularity of solutions to \eqref{eq:Homo-peri} over the space of periodic distributions, we consider the identification $\mathbf{U}(t)=\mathbf{u}(\cdot,t)$, with $\mathbf{U}:[0,\infty)\rightarrow \mathcal{H}^q(\T^n)$, for some $q\in\mathbb{R}$. In the following sections, spatial regularity refers to $\mathbf{U}(t)$ being in $\mathcal{H}^q(\T^n)$, for some $q\in\mathbb{R}$, for a fixed $t\ge 0$. Temporal regularity refers to the differentiability of the map $\mathbf{U}$ as a function of time $t$.

\begin{remark}
For any nonzero $k$, $M_k$ is a real symmetric invertible matrix. Hence, $M_k$ can be orthogonally diagonalized as
\begin{align*}
    M_k=P_k\begin{bmatrix}
    \lambda_1(k) & 0& \cdots&0\\
    0& \lambda_2(k) &\cdots &0\\
    \cdots&\cdots & \ddots & \cdots\\
    0& 0& \cdots& \lambda_2(k)
  \end{bmatrix}P_k^T,
\end{align*}
where $P_k$ is the matrix whose columns are the orthogonal eigenvectors of $M_k$. Then we can define the following matrices $\sqrt{-M_k},\; \cos(\sqrt{-M_k}t) $ and $\sin(\sqrt{-M_k}t)$ as
\begin{align*}
    \sqrt{-M_k}:=P_k\begin{bmatrix}
    \sqrt{-\lambda_1(k)} & 0& \cdots&0\\
    0& \sqrt{-\lambda_2(k)} &\cdots &0\\
    \cdots&\cdots & \ddots & \cdots\\
    0& 0& \cdots& \sqrt{-\lambda_2(k)}
  \end{bmatrix}P_k^T,
\end{align*}
and
\begin{align*}
    \cos(\sqrt{-M_k}t):=P_k\begin{bmatrix}
    \cos(\sqrt{-\lambda_1(k)}t) & 0& \cdots&0\\
    0& \cos(\sqrt{-\lambda_2(k)}t) &\cdots &0\\
    \cdots&\cdots & \ddots & \cdots\\
    0& 0& \cdots& \cos(\sqrt{-\lambda_2(k)}t)
  \end{bmatrix}P_k^T,
\end{align*}
and
\begin{align*}
    \sin(\sqrt{-M_k}t):=P_k\begin{bmatrix}
    \sin(\sqrt{-\lambda_1(k)}t) & 0& \cdots&0\\
    0& \sin(\sqrt{-\lambda_2(k)}t) &\cdots &0\\
    \cdots&\cdots & \ddots & \cdots\\
    0& 0& \cdots& \sin(\sqrt{-\lambda_2(k)}t)
  \end{bmatrix}P_k^T.
\end{align*}
From these definitions, one can show that for any $k\in \mathbb{Z}\setminus\{0\}$ and $t\ge 0$,
\begin{align*}
    \|(\sqrt{-M_k})^{-1}\|=\max\left\{\frac{1}{\sqrt{-\lambda_1(k)}},\frac{1}{\sqrt{-\lambda_2(k)}}\right\}=\frac{1}{\sqrt{\min\{|\lambda_1(k)|,|\lambda_2(k)|\}}},
\end{align*}
and
\begin{align*}
    \|\cos(\sqrt{-M_k}t)\|=\max\{|\cos(\sqrt{-\lambda_1(k)}t)|,|\cos(\sqrt{-\lambda_2(k)}t)|\},\\
    \|\sin(\sqrt{-M_k}t)\|=\max\{|\sin(\sqrt{-\lambda_1(k)}t)|,|\sin(\sqrt{-\lambda_2(k)}t)|\}.
\end{align*}
These observations imply that $\|\cos(\sqrt{-M_k}t)\|\le 1$ and $\|\sin(\sqrt{-M_k}t)\| \le 1$ for any $k\in \mathbb{Z}\setminus\{0\}$ and $t\ge 0$. Taking the classical derivative with respect to $t$, we have that
\begin{align*}
    &\frac{d^p}{d t^p}\begin{bmatrix}
    \cos(\sqrt{-\lambda_1(k)}) & 0& \cdots&0\\
    0& \cos(\sqrt{-\lambda_2(k)}) &\cdots &0\\
    \cdots&\cdots & \ddots & \cdots\\
    0& 0& \cdots& \cos(\sqrt{-\lambda_2(k)})
  \end{bmatrix}\\
  &=\begin{bmatrix}
    \cos(\sqrt{-\lambda_1(k)}+\frac{p\pi}{2}) & 0& \cdots&0\\
    0&\cos(\sqrt{-\lambda_2(k)}+\frac{p\pi}{2}) &\cdots &0\\
    \cdots&\cdots & \ddots & \cdots\\
    0& 0& \cdots& \cos(\sqrt{-\lambda_2(k)}+\frac{p\pi}{2})
  \end{bmatrix}(\sqrt{-M_k})^p.
\end{align*}
This implies that
\begin{equation*}
    \frac{d^p}{dt^p} \cos(\sqrt{-M_k}t)=\cos(\sqrt{-M_k}t+p\pi/2)(\sqrt{-M_k})^p.
\end{equation*}
Similarly, we  also have
\begin{align*}
    \frac{d^p}{d t^p}\sin(\sqrt{-M_k}t):=\sin(\sqrt{-M_k}t+p\pi/2)(\sqrt{-M_k})^p.
\end{align*}

\end{remark}

\subsection{Spatial and temporal regularity over periodic distributions}

Consider two fixed vector fields $\mathbf{f}\in \mathcal{H}^{s_1}(\mathbf{T}^n)$ and $\mathbf{g}\in \mathcal{H}^{s_2}(\mathbf{T}^n)$. Then for any $t\ge 0$, define
\begin{align}\label{Eq_wave}
    \mathbf{U}(t):=\sum_{k\in\mathbb{Z}^n}\widehat{\mathbf{U}}_k(t)e^{ik\cdot x}=
    (\widehat{\mathbf{f}}_0+\widehat{\mathbf{g}}_0 t)+\sum_{0\ne k\in \mathbb{Z}^n}\left(\cos(\sqrt{-M_k}t)\widehat{\mathbf{f}}_k+\sin(\sqrt{-M_k}t)(\sqrt{-M_k})^{-1}\widehat{\mathbf{g}}_k\right)e^{ik\cdot x}.
\end{align}
We will show in Theorem \ref{thm:Homo_existence} that \eqref{Eq_wave} is the unique solution to the peridynamic equation. We proceed by studying the spatial and temporal regularity of $\mathbf{U}$.

\begin{theorem}\label{thm:homo_1}
Let $n\ge 1$, $\delta>0$ and $\beta<n+2$. For any $t\ge 0$, $\mathbf{U}(t)\in \mathcal{H}^s(\mathbf{T}^n)$ where $s=\min\{s_1,s_2+\theta\}$ and $\theta=\max\{0,\frac{\beta-n}{2}\}$.
\end{theorem}

\begin{proof}
Notice that
\begin{align*}
    \sum_{k\in\mathbb{Z}^n}(1+\|k\|^2)^s|\hat{\mathbf{U}}_k(t)|^2\le &\sum_{0\ne k\in \mathbb{Z}^n}(1+\|k\|^2)^s|\widehat{\mathbf{f}}_k|^2\\
    +&\sum_{0\ne k\in \mathbb{Z}^n}(1+\|k\|^2)^s\|(\sqrt{-M_k})^{-1}\|^2\|\widehat{\mathbf{g}}_k\|^2+\|\widehat{\mathbf{f}_0}\|^2+\|\widehat{\mathbf{g}}_0\|^2t^2.
\end{align*}
Due to $f\in \mathcal{H}^{s_1}$ and $s\le s_1$, we get
\begin{align*}
    \sum_{0\ne k\in \mathbb{Z}^n}(1+\|k\|^2)^s\|\widehat{\mathbf{f}}_k\|^2\le \sum_{0\ne k\in \mathbb{Z}^n}(1+\|k\|^2)^{s_1}\|\widehat{\mathbf{f}}_k\|^2<\infty.
\end{align*}
Similarly, observe that
\begin{align*}
    \sum_{0\ne k\in \mathbb{Z}^n}(1+\|k\|^2)^s\|(\sqrt{-M_k})^{-1}\|^2\|\widehat{\mathbf{g}}_k\|^2=\sum_{0\ne k\in \mathbb{Z}^n}(1+\|k\|^2)^{s-s_2}\|(\sqrt{-M_k})^{-1}\|^2(1+\|k\|^2)^{s_2}|\hat{\mathbf{g}}_k|^2.
\end{align*}
Hence, the proof will be complete provided that
\begin{align*}
    (1+\|k\|^2)^{s-s_2}\|(\sqrt{-M_k})^{-1}\|^2,
\end{align*}
is bounded for large values of $\|k\|$.\\
To show this, we will consider two cases. First, when $\beta\le n$, then $s\le s_2$ and by Theorem~\ref{thm:asymptotics_2} and Theorem~\ref{thm:asymptotics_1}, there exists $r_1>0$ and $C_1>0$ such that   $\|(\sqrt{-M_k})^{-1}\|^2\le C_1$ for $\|k\|\ge r_1$. This leads to
\begin{align*}
    (1+\|k\|^2)^{s-s_2}\|(\sqrt{-M_k})^{-1}\|^2\le C_1.
\end{align*}
Similarly, for $\beta> n$, then $s \le s_2+\frac{\beta-n}{2}$ and by Theorem~\ref{thm:asymptotics_2} and Theorem~\ref{thm:asymptotics_1}, there exists $r_2>0$ and $C_2>0$ such that  $\|(\sqrt{-M_k})^{-1}\|^2\le C_2/( \|k\|^{\beta-n})$, for $\|k\|\ge r_2$. This leads to
\begin{align*}
    (1+\|k\|^2)^{s-s_2}\|(\sqrt{-M_k})^{-1}\|^2\le \frac{1}{C_2}\left(\frac{1+\|k\|^2}{\|k\|^2}\right)^{\frac{\beta-n}{2}},
\end{align*}
which is bounded.
\end{proof}

\begin{theorem}\label{thm:homo_reg}
Let $n\ge 1$, $\delta>0$ and $\beta<n+2$. Let $\mathbf{U}(t)$ be the map given by \eqref{Eq_wave}. Thus,
\begin{enumerate}
    \item if $\beta<n$, then $\mathbf{U}(t)\in C^\infty([0,\infty),\mathcal{H}^q(\mathbf{T}^n))$, for any $q\le s$,
    \item if $\beta=n$, then $\mathbf{U}(t)\in C^\infty([0,\infty),\mathcal{H}^q(\mathbf{T}^n))$, for any $q< s$, and
    \item if $\beta>n $, then $\mathbf{U}(t)\in C^{p+1}([0,\infty),\mathcal{H}^q(\mathbf{T}^n))$, for any $q\in \mathbb{R}$ and any positive integer $p$ satisfying 
    \begin{align*}
        q-s_1+(p+1)\frac{\beta-n}{2}\le 0\quad\text{and}\quad q-s_2+p\frac{\beta-n}{2}\le 0.
    \end{align*}
\end{enumerate}
\end{theorem}

\begin{proof}
We will prove this result by induction. Suppose $\mathbf{U}$ is already differentiable up to $p$ times, then by Lemma \ref{lem:derivative}, we have
\begin{align*}
    \mathbf{U}^{(p)}(t)=\sum_{k\in\mathbb{Z}^n}\widehat{\mathbf{U}}^{(p)}_k(t)e^{ik\cdot x},
\end{align*}
where $\widehat{\mathbf{U}}^{(p)}_k(t)$ is given by
\begin{align*}
    \widehat{\mathbf{U}}^{(p)}_k(t)=\cos\left(\sqrt{-M_k}t+\frac{p\pi}{2}\right)(\sqrt{-M_k})^p\widehat{\mathbf{f}}_k+\sin\left(\sqrt{-M_k}t+\frac{p\pi}{2}\right)(\sqrt{-M_k})^{p-1}\widehat{\mathbf{g}}_k,\quad \text{ for } k\ne 0
\end{align*}
and $\hat{U}^{(p)}_0(t)=\widehat{\mathbf{g}}_0\frac{d^p}{dt^p}(t)$. We will show that 
\begin{align*}
    \mathbf{U}^{(p+1)}(t)=\sum_{k\in\mathbb{Z}^n}\widehat{\mathbf{U}}^{(p+1)}_k(t)e^{ik\cdot x}
    .
\end{align*}
Using the mean value theorem, we have that
\begin{align*}
    T_h:&=\left\|\frac{\mathbf{U}^{(p)}(t+h)-\mathbf{U}^{(p)}(t)}{h}-\sum_{k\in\mathbb{Z}^n}\widehat{\mathbf{U}}^{(p+1)}_k(t)e^{ik\cdot x}\right\|^2_{\mathcal{H}^q(\T^n)}\\
    &=\sum_{k\in\mathbb{Z}^n}(1+\|k\|^2)^q\left|\frac{\widehat{\mathbf{U}}^{(p)}_k(t+h)-\widehat{\mathbf{U}}^{(p)}_k(t)}{h}-\widehat{\mathbf{U}}^{p+1}_k(t)\right|^2\\
    &\le\sum_{0\ne k\in\mathbb{Z}^n}(1+\|k\|^2)^q\left\|\cos\left(\sqrt{-M_k}(t+\xi_{p,k})+\frac{(p+1)\pi}{2}\right)\right.\\
    &\qquad \qquad \qquad \left. -\cos\left(\sqrt{-M_k}t+\frac{(p+1)\pi}{2}\right)\right\|^2\|\sqrt{-M_k}\|^{2(p+1)}\|\widehat{\mathbf{f}}_k\|^2\\
    &+\sum_{0\ne k\in\mathbb{Z}^n}(1+\|k\|^2)^q|\left\|\sin\left(\sqrt{-M_k}(t+\eta_{p,k})+\frac{(p+1)\pi}{2}\right)\right.\\
    &\qquad \qquad \qquad\left.-\sin\left(\sqrt{-M_k}t+\frac{(p+1)\pi}{2}\right)\right\|^2\|\sqrt{-M_k}\|^{2p}\|\widehat{\mathbf{g}}_k\|^2
\end{align*}
for some $\xi_{p,k},\eta_{p,k}\in(0,h)$. This means that as $h\to 0$, each term of the summation will approach 0. The proof will be complete if we are able to pass the limit inside the summation. This is achievable provided the uniform boundedness of
\begin{align*}
    a_k:=(1+\|k\|^2)^{q-s_1}\|\sqrt{-M_k}\|^{2(p+1)}\quad\text{and}\quad b_k:=(1+\|k\|^2)^{q-s_2}\|\sqrt{-M_k}\|^{2p}
\end{align*}
To establish these bounds, we consider the three cases when $\beta<n$, $\beta=n$ and $\beta>n$ respectively.\\
When $\beta<n$, from Theorem~\ref{thm:asymptotics_2} and Theorem`\ref{thm:asymptotics_1}, $\|\sqrt{-M_k}\|$ is bounded, and this  provides the uniform boundedness of $a_k$ and $b_k$ provided that $q\le s_1,s_2$.\\
When $\beta=n$, fix an $\epsilon>0$. Theorem~\ref{thm:asymptotics_2} and Theorem`\ref{thm:asymptotics_1} show that there exists constants $C,r>0$ such that $\|\sqrt{-M_k}\|^2\le C\log\|k\|\le C(1+\|k\|^2)^{\epsilon/2}$ whenever $\|k\|\ge r$. Then for large $\|k\|$,
\begin{align*}
    a_k\le C(1+\|k\|^2)^{q-s_1+(p+1)\epsilon/2}\quad\text{and}\quad b_k\le C(1+\|k\|^2)^{q-s_2+p\epsilon/2}.
\end{align*}
This offers the uniform boundedness of $a_k$ and $b_k$ provided that $q-s_1+(p+1)\epsilon/2\le 0$ and $q-s_2+p\epsilon/2$. So given $q<s_1,s_2$, we can always pick $\epsilon$ satisfying these conditions.\\
When $\beta>n$, again from Theorem~\ref{thm:asymptotics_2} and Theorem`\ref{thm:asymptotics_1}, there is $C,r>0$ such that $\|\sqrt{-M_k}\|^2\le C\|k\|^{\beta-n}\le C(1+\|k\|^2)^{(\beta-n)/2}$ for $\|k\|\ge r$. Then for large $\|k\|$
\begin{align*}
    a_k\le C(1+\|k\|^2)^{q-s_1+(p+1)\frac{\beta-n}{2}}\quad\text{and}\quad b_k\le C(1+\|k\|^2)^{q-s_2+p\frac{\beta-n}{2}}
\end{align*}
This shows that the uniform boundedness of $a_k,b_k$ is established given the conditions
\begin{align*}
    q-s_1+(p+1)\frac{\beta-n}{2}\le 0\quad\text{and}\quad q-s_2+p\frac{\beta-n}{2}\le 0.
\end{align*}
\end{proof}

Theorem~\ref{thm:Homo_existence} is a direct application of Theorem~\ref{thm:homo_1} and Theorem~\ref{thm:homo_reg}. The assumptions are to guarantee that the second derivative exists.
\begin{theorem}\label{thm:Homo_existence}
Let $n\ge 1$, $\delta>0$ and $\beta<n+2$. Let $\mathbf{f}\in \mathcal{H}^{s_1}(\T^n)$ and $\mathbf{g}\in \mathcal{H}^{s_2}(\T^n)$. Suppose that
\begin{enumerate}
    \item $q\le \min\{s_1,s_2\}$, in the case $\beta<n$,
    \item $q<\min\{s_1,s_2\}$, in the case $\beta=n$, and
    \item $q-s_1+(\beta-n)\le 0$ and $q-s_2+(\beta-n)/2\le 0$, in the case $\beta>n$.
\end{enumerate}
Then the map $\mathbf{U}(t)$ defined by \eqref{Eq_wave} is in $\mathcal{H}^q(\mathbf{T}^n)$, for $t\ge 0$. Moreover, $\mathbf{U}$ is the unique solution to the peridynamic equation
\begin{align}\label{eq:Homo-peri-new}
    \begin{cases}
        \frac{d^2}{dt^2}\mathbf{U}(t)&=\Lper^{\delta,\beta}\mathbf{U}(t),\; t>0\\
    \mathbf{U}(0)&=\mathbf{f},\\
    \frac{d}{dt}\mathbf{U}(0)&=\mathbf{g}.
    \end{cases}
\end{align}
Furthermore, when $\beta\le n$, $\mathbf{U}\in C^\infty([0,\infty); \mathcal{H}^s(\T^n))$ where $s=\min\{s_1,s_2\}$,  and when \\$\beta>n$, 
$\mathbf{U}\in C^2([0,\infty); \mathcal{H}^s(\T^n))$ where $s=\min\{s_1,s_2+\frac{\beta-n}{2}\}$.

\end{theorem}

\subsection{Convergence to the transient Navier-Cauchy equation}\label{sec:Conv-homo}
In this section, we will provide some results on the convergence of the solutions of the peridynamic equation \eqref{eq:Homo-peri-new} to the  solution of the corresponding Navier-Cauchy equation for two types of  limits as $\delta \to 0^+$, with $\beta< n+2$ is being fixed,  or as $\beta\to (n+2)^-$, with $\delta>0$ is being fixed. The transient Navier-Cauchy equation is given by
\begin{align}\label{eq:Homo-Navier}
    \begin{cases}
        \frac{d^2}{dt^2}\mathbf{U}^\mathcal{N}(t)&=\mathcal{N}\mathbf{U}^{\mathcal{N}}(t),\; t>0\\
    \mathbf{U}^{\mathcal{N}}(0)&=\mathbf{f},\\
    \frac{d}{dt}\mathbf{U}^{\mathcal{N}}(0)&=\mathbf{g},
    \end{cases}
\end{align}
and its solution is given by
\begin{align}
    \mathbf{U}^{\mathcal{N}}(t):=
    (\widehat{\mathbf{f}}_0+\widehat{\mathbf{g}}_0 t)+\sum_{0\ne k\in \mathbb{Z}^n}\left(\cos(\sqrt{-M^{\mathcal{N}}_k}t)\widehat{\mathbf{f}}_k+\sin(\sqrt{-M^{\mathcal{N}}_k}t)(\sqrt{-M^{\mathcal{N}}_k})^{-1}\widehat{\mathbf{g}}_k\right)e^{ik\cdot x}.
\end{align}
To emphasize the dependence of the solution of the nonlocal wave equation \eqref{eq:Homo-peri-new} on the parameters $\delta$ and $\beta$, in the remaining part of this section we denote the solution $\mathbf{U}$ by $\mathbf{U}^{\delta,\beta}$.


\begin{theorem}\label{thm:homo1}
Let $n\ge 1$, $\delta>0$ and $\beta<n+2$. Suppose $\mathbf{f}\in \mathcal{H}^{s_1}(\T^n)$ and $\mathbf{g}\in \mathcal{H}^{s_2}(\T^n)$. Let $\mathbf{U}^{\delta,\beta}(t)\in \mathcal{H}^{s'}(\T^n)$ and $\mathbf{U}^{\mathcal{N}}(t)\in \mathcal{H}^s(\T^n)$  be the solutions to the peridynamic equation \eqref{eq:Homo-peri-new} and Navier-Cauchy equation \eqref{eq:Homo-Navier}, respectively, with $s'=\min\{s_1,s_2+\max\{0,\frac{\beta-n}{2}\}\}$ and $s=\min\{s_1,s_2+1\}$. Let $t\ge 0$, then when $\beta<n+2$ is fixed, 
\begin{align*}
    \lim_{\delta\to 0^+}\mathbf{U}^{\delta,\beta}(t)=\mathbf{U}^{\mathcal{N}}(t),\quad \text{in }\mathcal{H}^{s'}(\T^n).
\end{align*}
\end{theorem}

\begin{proof}
From \ref{Eq_wave}, one can show that
\begin{align*}
    \|\mathbf{U}^{\delta,\beta}(t)-&\mathbf{U}^{\mathcal{N}}(t)\|^2_{\mathcal{H}^{s'}(\T^n)}\le \sum_{k\in \mathbb{Z}^n}(1+\|k\|^2)^{s'}\left\|\cos(\sqrt{-M_k}t)-\cos(\sqrt{-M_k^{\mathcal{N}}}t)\right\|^2\|\widehat{\mathbf{f}}_k\|^2\\
    &+\sum_{k\in \mathbb{Z}^n}(1+\|k\|^2)^{s'}\left\|\sin(\sqrt{-M_k}t)(\sqrt{-M_k})^{-1}-\sin(\sqrt{-M_k^{\mathcal{N}}}t)(\sqrt{-M_k^{\mathcal{N}}})^{-1}\right\|^2\|\widehat{\mathbf{g}}_k\|^2
\end{align*}
Using  Proposition \ref{prop1}, each term tends to $0$ under the limits, hence, the proof is complete if we can pass the limit inside the summation. Since $\mathbf{f}\in \mathcal{H}^{s_1}(\T^n)$ and $\mathbf{g}\in \mathcal{H}^{s_2}(\T^n)$, this is achievable by establishing the uniform boundedness for 
\[
    \|k\|^{s'-s_1}\left\|\cos(\sqrt{-M_k}t)-\cos(\sqrt{-M_k^{\mathcal{N}}}t)\right\|
\]
and
\[
    \|k\|^{s'-s_2}\left\|\sin(\sqrt{-M_k}t)(\sqrt{-M_k})^{-1}-\sin(\sqrt{-M_k^{\mathcal{N}}}t)(\sqrt{-M_k^{\mathcal{N}}})^{-1}\right\|
\]
as $\|k\|\to \infty$. \\
It is straightforward to see that
\begin{align*}
    \left\|\cos(\sqrt{-M_k}t)-\cos(\sqrt{-M_k^{\mathcal{N}}}t)\right\|=\max\limits_{i=1,2}\left|\cos(\sqrt{-\lambda_i(k)}t)-\cos(\sqrt{-\lambda_i^{\mathcal{N}}(k)}t)\right|\le 2.
\end{align*}
This gives the bound for the first term since $s'\le s_1$. From Theorem \ref{thm:asymptotics_2} and Theorem \ref{thm:asymptotics_1}, there exists $r,C_1,C_2>0$ such that
\begin{align*}
    |\lambda_1(k)|,|\lambda_2(k)|\ge\begin{cases}
    C_1&\text{if }\beta\le n\\
    C_2\|k\|^{\beta-n}&\text{if }\beta>n
    \end{cases}
\end{align*}
for any $\|k\|\ge r$. Note that
\begin{align*}
    &\left\|\sin(\sqrt{-M_k}t)(\sqrt{-M_k})^{-1}-\sin(\sqrt{-M_k^{\mathcal{N}}}t)(\sqrt{-M_k^{\mathcal{N}}})^{-1}\right\| \\
&=\max\limits_{i=1,2}\left|\frac{\sin(\sqrt{-\lambda_i(k)}t)}{\sqrt{-\lambda_i(k)}}-\frac{\sin(\sqrt{-\lambda_i^{\mathcal{N}}(k)}t)}{\sqrt{-\lambda^{\mathcal{N}}_i(k)}}\right|\\
    &\le\max\limits_{i=1,2}\left(\frac{1}{\sqrt{-\lambda_i(k)}}+\frac{1}{\sqrt{-\lambda^{\mathcal{N}}_i(k)}}\right)\\
    &\le \begin{cases}
    \frac{D}{\|k\|}+C&\text{if }\beta\le n\\
    \frac{D}{\|k\|}+\frac{C}{\|k\|^{\frac{\beta-n}{2}}}&\text{if }\beta> n
    \end{cases}.
\end{align*}
This result combined with the fact that
\begin{align*}
    \|k\|^{s'-s_2}\le\begin{cases}
    1&\text{if }\beta\le n\\
    \|k\|^{\frac{\beta-n}{2}}&\text{if }\beta>n
    \end{cases},
\end{align*}
yields the bound for the second term.
\end{proof}

\begin{theorem}
Let $n\ge 1$, $\delta>0$ and $\beta<n+2$. Suppose $\mathbf{f}\in \mathcal{H}^{s_1}(\T^n)$ and $\mathbf{g}\in \mathcal{H}^{s_2}(\T^n)$. Let $\mathbf{U}^{\delta,\beta}(t)\in \mathcal{H}^{s'}(\T^n)$ and $\mathbf{U}^{\mathcal{N}}(t)\in \mathcal{H}^s(\T^n)$  be the solutions to the peridynamic equation \eqref{eq:Homo-peri-new} and Navier-Cauchy equation \eqref{eq:Homo-Navier}, respectively, with $s'=\min\{s_1,s_2+\max\{0,\frac{\beta-n}{2}\}\}$ and $s=\min\{s_1,s_2+1\}$. Let $t\ge 0$, then, for $0<\epsilon<2$,
\begin{align*}
    \lim_{\beta\to n+2^-}\mathbf{U}^{\delta,\beta}(t)=\mathbf{U}^{\mathcal{N}}(t),\quad \text{in }\mathcal{H}^{s_0}(\T^n),
\end{align*}
where $s_0=\min\{s_1,s_2+\frac{2-\epsilon}{2}\}$.
\end{theorem}

\begin{proof}
Fix an $0<\epsilon<2$, for any $\beta<n+2$ such that $n+2-\beta<\epsilon$ both $\mathbf{U}(t)$ and $\mathbf{U}^{\mathcal{N}}(t)$ belong to the space $\mathcal{H}^{s_0}(\T^n)$. With the same approach as in the proof of Theorem \ref{thm:homo1}, we can show that
\begin{align*}
    \|\mathbf{U}^{\delta,\beta}(t)-&\mathbf{U}^{\mathcal{N}}(t)\|^2_{\mathcal{H}^{s_0}(\T^n)}\le \sum_{k\in \mathbb{Z}^n}(1+\|k\|^2)^{s_0}\left\|\cos(\sqrt{-M_k}t)-\cos(\sqrt{-M_k^{\mathcal{N}}}t)\right\|^2\|\widehat{\mathbf{f}}_k\|^2\\
    &+\sum_{k\in \mathbb{Z}^n}(1+\|k\|^2)^{s_0}\left\|\sin(\sqrt{-M_k}t)(\sqrt{-M_k})^{-1}-\sin(\sqrt{-M_k^{\mathcal{N}}}t)(\sqrt{-M_k^{\mathcal{N}}})^{-1}\right\|^2\|\widehat{\mathbf{g}}_k\|^2.
\end{align*}
Next we establish uniform bounds for the following terms
\begin{align*}
    \|k\|^{s_0-s_1}\left\|\cos(\sqrt{-M_k}t)-\cos(\sqrt{-M_k^{\mathcal{N}}}t)\right\|
\end{align*}
and
\begin{align*}
    \|k\|^{s_0-s_2}\left\|\sin(\sqrt{-M_k}t)(\sqrt{-M_k})^{-1}-\sin(\sqrt{-M_k^{\mathcal{N}}}t)(\sqrt{-M_k^{\mathcal{N}}})^{-1}\right\|
\end{align*}
The bound for the first term is the same as in Theorem \ref{thm:homo1}. For the second term, since $\beta>n$, we have shown that
\begin{align*}
    \left\|\sin(\sqrt{-M_k}t)(\sqrt{-M_k})^{-1}-\sin(\sqrt{-M_k^{\mathcal{N}}}t)(\sqrt{-M_k^{\mathcal{N}}})^{-1}\right\|\le 
    \frac{D}{\|k\|}+\frac{C}{\|k\|^{\frac{\beta-n}{2}}}.
\end{align*}
Since $\|k\|^{s_0-s_2}\le\|k\|^{\frac{2-\epsilon}{2}}$, this implies that
\begin{align*}
    \|k\|^{s_0-s_2}\left\|\sin(\sqrt{-M_k}t)(\sqrt{-M_k})^{-1}-\sin(\sqrt{-M_k^{\mathcal{N}}}t)(\sqrt{-M_k^{\mathcal{N}}})^{-1}\right\|\le \frac{D}{\|k\|^{\frac{\epsilon}{2}}}+\frac{C}{\|k\|^{\frac{\epsilon+\beta-(n+2)}{2}}},
\end{align*}
which is bounded.
\end{proof}

\subsection{Spatial and temporal regularity over \texorpdfstring{$(L^2(\T^n))^n$}{}}\label{sec:temp-homo}
The solution $\mathbf{U}(t)\in \mathcal{H}^s(\T^n)$ given in \eqref{Eq_wave} is a distribution when $s<0$, and  when $s\ge 0$ defines a regular vector field
\begin{align*}
     \mathbf{u}(x,t)=\mathbf{U}(t)(x)=
    (\widehat{\mathbf{f}}_0+\widehat{\mathbf{g}}_0 t)+\sum_{0\ne k\in \mathbb{Z}^n}\left(\cos(\sqrt{-M_k}t)\widehat{\mathbf{f}}_k+\sin(\sqrt{-M_k}t)(\sqrt{-M_k})^{-1}\widehat{\mathbf{g}}_k\right)e^{ik\cdot x}.
\end{align*}
Let the data be two regular vector fields: $\mathbf{f}\in \mathcal{H}^{s_1}(\T^n)$ and $\mathbf{g}\in \mathcal{H}^{s_2}(\T^n)$ with $s_1\ge 0$ and $s_2\ge 0$. In this section, we focus on conditions that guarantee that $\mathbf{u}(x,t)$, the solution of \eqref{eq:Homo-peri}, is a regular vector field $\mathbf{u}(\cdot,t)\in (L^2(\T^n))^n$. 

\subsubsection{Functions with absolutely summable Fourier coefficients}
In this section, we state some well-known results on the differentiablity and summability of series in multidimensions. These results will be used in subsequent sections and are included here for completeness of the presentation. 

We recall that a vector field  $\mathbf{f}$ on $\mathbb{R}^n$ is called $\alpha$-H{\"o}lder continuous, where $\alpha \in [0,1]$, if there exists a constant $C$ such that 
\begin{align*}
    \| \mathbf{f}(x)-\mathbf{f}(y)| \|\le C\|x-y\|^\alpha,
\end{align*}
for any $x,y$ in the domain of $\mathbf{f}$.

\begin{lemma}\label{Fourier_summable}
Suppose that $\mathbf{f}$ is a vector field defined on the torus $\T^n$ satisfying the $\alpha$-H{\"o}lder continuity condition for some $\alpha>n/2$. Then 
\begin{align*}
\sum_{k\in\mathbb{Z}^n}\|\widehat{\mathbf{f}}_k\|<\infty.
\end{align*}
\end{lemma}
A version of the above result for scalar functions can be found in \cite{Grafakos2009}.

\begin{lemma}\label{lem:series}
Let $T$ be an open subset of $\mathbb{R}^+$, $K$ is a measurable space equipped with the counting measure and $\mathbf{f}_k:T\to \mathbb{C}^n$. Suppose that
\begin{enumerate}
    \item $\sum_{k\in K}\|\mathbf{f}_k(t)\|<\infty$ for each $t\in T$.
    \item For almost all $k\in K$, the derivative $\frac{d}{dt}\mathbf{f}_k(t)$ exists and is continuous for all $t\in T$.
    \item There is a summable sequence $\theta_k$ such that $\|\frac{d}{dt}\mathbf{f}_k(t)\|\le \theta_k$ for all $t\in T$ and almost all $k\in K$.
\end{enumerate}
Then, for all $t\in T$,
\begin{align*}
    \frac{d}{dt}\sum_{k}\mathbf{f}_k(t)=\sum_{k}\frac{d}{dt}\mathbf{f}_k(t).
\end{align*}
\end{lemma}
A similar statement to Lemma \ref{lem:series} for scalar functions can be found in \cite{Folland1999}.

\subsubsection{Temporal regularity with respect to the Gateaux derivative}
Let the data $\mathbf{f}\in \mathcal{H}^{s_1}(\T^n)$ and $\mathbf{g}\in \mathcal{H}^{s_2}(\T^n)$, with $s_1\ge 0$ and $s_2\ge 0$.
\begin{theorem}
Let $n\ge 1$, $\delta>0$, and $\beta<n+2$. 
Suppose that
\begin{enumerate}
    \item $s_1,s_2\ge 0$, in the case when $\beta<n$,
    \item $s_1,s_2>0$, in the case when $\beta=n$, and
    \item $s_1\ge \frac{3}{2}(\beta-n) $ and $s_2\ge \beta-n$, in the case when $\beta>n$.
\end{enumerate}
Then, the  vector field $\mathbf{u}$, where $  \mathbf{u}(x,t)=\mathbf{U}(t)(x)$, with $\mathbf{U}$ defined  in \eqref{Eq_wave}, is the unique solution to the  peridynamic equation \eqref{eq:Homo-peri}.
Moreover, when $\beta\le n$, then $\mathbf{u}(x,t)\in C^\infty([0,\infty); \mathcal{H}^s(\T^n))$, where $s=\min\{s_1,s_2\}$, and when $\beta>n$, 
then $\mathbf{u}(x,t)\in C^2([0,\infty); \mathcal{H}^s(\T^n))$, where
$s=\min\{s_1,s_2+\frac{\beta-n}{2}\}$.
\end{theorem}

\subsubsection{Temporal regularity with respect to the classical derivative}

\begin{theorem}
    Let $n\ge 1$, $\delta>0$ and $\beta<n$. Suppose that $\mathbf{f}$ and $\mathbf{g}$ satisfy the $\alpha$ - H{\"o}lder continuity condition with $\alpha>n/2$. Then $\mathbf{u}(x,\cdot)\in C^\infty[0,\infty)$, for any $x\in\T^n$.
\end{theorem}

\begin{proof}
Fix a value of $p$ in $\mathbb{Z}_{\ge 0}$. Consider the two series 
\begin{align*}
    \displaystyle \sum_{0\ne k\in \mathbb{Z}^n}\cos(\sqrt{-M_k}t)\widehat{\mathbf{f}}_k e^{ik\cdot x} \quad\text{and}\quad \displaystyle \sum_{0\ne k\in \mathbb{Z}^n}\sin(\sqrt{-M_k}t)(\sqrt{-M_k})^{-1}\widehat{\mathbf{g}}_k e^{ik\cdot x}.
\end{align*}
From Theorem \ref{thm:asymptotics_2} and Theorem \ref{thm:asymptotics_1}, $\|\sqrt{-M_k}\|$ and $\|(\sqrt{-M_k})^{-1}\|$ are bounded for large $k$. This makes
\begin{align*}
    &\|\cos(\sqrt{-M_k}t)\widehat{\mathbf{f}}_k e^{ik\cdot x}\|\le\|\widehat{\mathbf{f}}_k\|,\\
    &\left\|\frac{d^p}{dt^p}\cos(\sqrt{-M_k}t)\widehat{\mathbf{f}}_k e^{ik\cdot x}\right\|\le\|\sqrt{-M_k}\|^p\|\widehat{\mathbf{f}}_k\|\le C_1  \|\widehat{\mathbf{f}}_k\|,\\
    &\|\sin(\sqrt{-M_k}t)(\sqrt{-M_k})^{-1}\widehat{\mathbf{g}}_ke^{ik\cdot x}\|\le\|(\sqrt{-M_k})^{-1}\|\|\widehat{\mathbf{g}}_k\|\le C_2 \|\widehat{\mathbf{g}}_k\|,\\
    &\left\|\frac{d^p}{dt^p}\sin(\sqrt{-M_k}t)(\sqrt{-M_k})^{-1}\widehat{\mathbf{g}}_ke^{ik\cdot x}\right\|\le\|(\sqrt{-M_k})\|^{p-1}\|\widehat{\mathbf{g}}_k\|\le C_3 \|\widehat{\mathbf{g}}_k\|,
\end{align*}
as $k\to \infty$. From Lemma \ref{Fourier_summable}, the Fourier coefficients of $\widehat{\mathbf{f}}_k$ and $\widehat{\mathbf{g}}_k$ are absolutely summable. The proof is now complete by using Lemma \ref {lem:series}.
\end{proof}

\begin{theorem}
    Let $n\ge 1$, $\delta>0$ and $\beta=n$. Suppose $s_1>n$ and $s_2>n$. Then $\mathbf{u}(x,\cdot)\in C^\infty[0,\infty)$, for any $x\in\T^n$. 
\end{theorem}

\begin{proof}
    Fix a value of $p$ in $\mathbb{Z}_{\ge 0}$ and choose $\epsilon>0$ such that $s_1-p\epsilon>n$ and $s_2-(p-1)\epsilon>n$. Consider the two series
    \begin{align*}
    \displaystyle \sum_{0\ne k\in \mathbb{Z}^n}\cos(\sqrt{-M_k}t)\widehat{\mathbf{f}}_k e^{ik\cdot x} \quad\text{and}\quad \displaystyle \sum_{0\ne k\in \mathbb{Z}^n}\sin(\sqrt{-M_k}t)(\sqrt{-M_k})^{-1}\widehat{\mathbf{g}}_k e^{ik\cdot x}.
\end{align*}
From Theorem \ref{thm:asymptotics_2} and Theorem \ref{thm:asymptotics_1}, 
\begin{align*}
    \|\sqrt{-M_k}\|\le C(\log\|k\|)^{1/2}\le C\|k\|^{\epsilon},\quad\text{and}\quad
    \|(\sqrt{-M_k})^{-1}\|\le  C,
\end{align*}
as $k\to \infty$. This implies 
\begin{align*}
    &\|\cos(\sqrt{-M_k}t)\widehat{\mathbf{f}}_k e^{ik\cdot x}\|\le\|\widehat{\mathbf{f}}_k\|,\\
    &\left\|\frac{d^p}{dt^p}\cos(\sqrt{-M_k}t)\widehat{\mathbf{f}}_k e^{ik\cdot x}\right\|\le\|\sqrt{-M_k}\|^p\|\widehat{\mathbf{f}}_k\|\le \frac{C_1}{\|k\|^{s_1-p\epsilon}},\\
    &\|\sin(\sqrt{-M_k}t)(\sqrt{-M_k})^{-1}\widehat{\mathbf{g}}_ke^{ik\cdot x}\|\le\|(\sqrt{-M_k})^{-1}\|\|\widehat{\mathbf{g}}_k\|\le \frac{C_2}{\|k\|^{s_2}} ,\\
    &\left\|\frac{d^p}{dt^p}\sin(\sqrt{-M_k}t)(\sqrt{-M_k})^{-1}\widehat{\mathbf{g}}_ke^{ik\cdot x}\right\|\le\|(\sqrt{-M_k})\|^{p-1}\|\widehat{\mathbf{g}}_k\|\le \frac{C_3}{\|k\|^{s_2-(p-1)\epsilon}},
\end{align*}
as $k\to \infty$. The proof is  complete by using Lemma \ref {lem:series}.
\end{proof}

\begin{theorem}
    Let $n\ge 1$, $\delta>0$ and $n<\beta<n+2$. Suppose that there exist a nonnegative integer $p$ such that 
    \begin{align*}
        s_1>n+p\frac{\beta-n}{2}\quad\text{and}\quad s_2>n+(p-1)\frac{\beta-n}{2}.
    \end{align*} Then $\mathbf{u}(x,\cdot)\in C^p[0,\infty)$, for any $x\in\T^n$.
\end{theorem}

\begin{proof}
     Consider the two series \begin{align*}
    \displaystyle \sum_{0\ne k\in \mathbb{Z}^n}\cos(\sqrt{-M_k}t)\widehat{\mathbf{f}}_k e^{ik\cdot x} \quad\text{and}\quad \displaystyle \sum_{0\ne k\in \mathbb{Z}^n}\sin(\sqrt{-M_k}t)(\sqrt{-M_k})^{-1}\widehat{\mathbf{g}}_k e^{ik\cdot x}.
 \end{align*} 
From Theorem \ref{thm:asymptotics_2} and Theorem \ref{thm:asymptotics_1}, 
\begin{align*}
    \|\sqrt{-M_k}\|\le C\|k\|^{(\beta-n)/2},\quad\text{and}\quad
    \|(\sqrt{-M_k})^{-1}\|\le  C\|k\|^{-(\beta-n)/2},
\end{align*}
as $k\to \infty$. This implies 
\begin{align*}
    &\|\cos(\sqrt{-M_k}t)\widehat{\mathbf{f}}_k e^{ik\cdot x}\|\le\|\widehat{\mathbf{f}}_k\|,\\
    &\left\|\frac{d^p}{dt^p}\cos(\sqrt{-M_k}t)\widehat{\mathbf{f}}_k e^{ik\cdot x}\right\|\le\|\sqrt{-M_k}\|^p\|\widehat{\mathbf{f}}_k\|\le \frac{C_1}{\|k\|^{s_1-p(\beta-n)/2}},\\
    &\|\sin(\sqrt{-M_k}t)(\sqrt{-M_k})^{-1}\widehat{\mathbf{g}}_ke^{ik\cdot x}\|\le\|(\sqrt{-M_k})^{-1}\|\|\widehat{\mathbf{g}}_k\|\le \frac{C_2}{\|k\|^{s_2+(\beta-n)/2}},\\
    &\left\|\frac{d^p}{dt^p}\sin(\sqrt{-M_k}t)(\sqrt{-M_k})^{-1}\widehat{\mathbf{g}}_ke^{ik\cdot x}\right\|\le\|(\sqrt{-M_k})\|^{p-1}\|\widehat{\mathbf{g}}_k\|\le \frac{C_3}{\|k\|^{s_2-(p-1)(\beta-n)/2}},
\end{align*}
as $k\to \infty$. The proof is  complete by using Lemma \ref {lem:series}.
\end{proof}

\section{The peridynamic equation with a forcing term }\label{sec:Peri-forcing}
In this section, we focus on the following nonlocal wave equation with a  forcing term and zero initial data. 
\begin{align}\label{eq:Nonhomo-peri}
    \begin{cases}
    \mathbf{u}_{tt}(x,t)&=\Lper \mathbf{u}(x,t)+\mathbf{b}(x),\quad x\in\T^n, \; t>0,\\
    \mathbf{u}(x,0)&=0,\\
    \mathbf{u}_t(x,0)&=0.
    \end{cases}
\end{align}
In order to study the existence, uniqueness, and regularity of solutions to \eqref{eq:Nonhomo-peri} over the space of periodic distributions, we consider the identification $\mathbf{U}(t)=\mathbf{u}(\cdot,t)$, with $\mathbf{U}:[0,\infty)\rightarrow H^s(\T^n)$.

\subsection{Spatial and temporal regularity over periodic distributions}
For a fixed vector field $\mathbf{b}\in \mathcal{H}^S(\T^n)$ and any $t\ge 0$, define
\begin{align}\label{eq:forcing}
    \mathbf{U}(t):=\sum_{k\in \mathbb{Z}^n}\widehat{\mathbf{U}}_k(t)e^{ik\cdot x}=\widehat{\mathbf{b}}_0\frac{t^2}{2}+\sum_{0\ne k\in \mathbb{Z}^n}\left[\cos\left(\sqrt{-M_k}t\right)-I\right]M_k^{-1}\widehat{\mathbf{b}}_k e^{ik\cdot x}.
\end{align}

\begin{theorem}\label{thm:for_1}
Let $n\ge 1$, $\delta>0$, and $\beta<n+2$. For any $t\ge 0$, $\mathbf{U}(t)\in \mathcal{H}^s(\T^n)$ where $s=S+\max\{0,\beta-n\}$.
\end{theorem}
\begin{proof}
We can show that
\begin{align*}
    \|\cos\left(\sqrt{-M_k}t\right)-I\|=\max\{|\cos(\sqrt{-\lambda_1(k)}t)-1|,|\cos(\sqrt{-\lambda_2(k)}t)-1|\}\le 2
\end{align*}
Then it follows that
\begin{align*}
    \sum_{k\in\mathbb{Z}^n}(1+\|k\|^2)^s\|\widehat{\mathbf{U}}_k(t)\|^2\le \|\hat{\mathbf{b}}_0\|\frac{t^2}{2}+\sum_{0\ne k\in \mathbb{Z}^n}4(1+\|k\|^2)^s\|M_k^{-1}\|^2\|\widehat{\mathbf{b}}_k\|^2
\end{align*}
Since $\mathbf{b}\in H^S(\T^n)$, the proof is complete provided the uniform boundedness of
\begin{align*}
    (1+\|k\|^2)^\theta\|M_k^{-1}\|^2,
\end{align*}
for large value of $\|k\|$, where $\theta:=\max\{0,\beta-n\}$. To show this, we consider the two cases $\beta\le n$ and $\beta > n$ respectively.\\
For the case $\beta\le n$, then $\theta=0$. From Theorem~\ref{thm:asymptotics_2} and Theorem~\ref{thm:asymptotics_1}, there are constants $C_1,r_1>0$ such that $\|M_k^{-1}\|^2\le C_1$ whenever $\|k\|\ge r_1$, this implies that
\begin{align*}
    (1+\|k\|^2)^\theta\|M_k^{-1}\|^2\le C_1.
\end{align*}
When $\beta>n$, then $\theta=\beta-n$. Also from Theorem~\ref{thm:asymptotics_2} and Theorem~\ref{thm:asymptotics_1}, there exists $C_2,r_2>0$ satisfying $\|M_k^{-1}\|^2\le C_2/\|k\|^{2(\beta-n)}$ for any $\|k\|\ge r_2$. This yields
\begin{align*}
    (1+\|k\|^2)^\theta\|M_k^{-1}\|^2\le C_2\left(\frac{1+\|k\|^2}{\|k\|^2}\right)^{\beta-n}
\end{align*}
\end{proof}

\begin{theorem}\label{thm:for_regularity}
Let $n\ge 1$, $\delta>0$, and $\beta<n+2$. Let $\mathbf{U}(t)$ be the map given by \eqref{eq:forcing}. Thus
\begin{enumerate}
    \item if $\beta<n$, then $\mathbf{U}(t)\in C^\infty([0,\infty),\mathcal{H}^q(\T^n))$, for any $q\le S$, 
    \item if $\beta=n$, $\mathbf{U}(t)\in C^\infty([0,\infty),\mathcal{H}^q(\T^n))$, for any $q<S$ , and
    \item if $\beta>n $, then $\mathbf{U}(t)\in C^{p+1}([0,\infty),\mathcal{H}^q(\T^n))$. for any $q\in \mathbb{R}$ and any positive integer $p$ satisfying 
    \begin{align*}
        q-S+(p-1)\frac{\beta-n}{2}\le 0.    \end{align*}
    
\end{enumerate}
\end{theorem}
\begin{proof}
We will show this result by induction. Suppose $\mathbf{U}(t)$ is already differentiable up to $p$ times. According to Lemma \ref{lem:derivative}, we have that 
\begin{align*}
    \mathbf{U}^{(p)}(t)=\sum_{k\in\mathbb{Z}^n}\widehat{\mathbf{U}}^{(p)}_k(t)e^{ik\cdot x},
\end{align*}
where $\widehat{\mathbf{U}}^{(p)}_k(t)$ is given by
\begin{align*}
    \widehat{\mathbf{U}}^{(t)}_k(t)=\cos\left(\sqrt{-M_k}t+\frac{p\pi}{2}\right)(\sqrt{-M_k})^{p-2} \widehat{\mathbf{b}}_k,
\end{align*}
and $\widehat{\mathbf{U}}^{(p)}_0(t)=\frac{\widehat{\mathbf{b}}_0}{2}\frac{d^p}{dt^p}(t^2)$. Using induction, we need to show that
\begin{align*}
    \mathbf{U}^{(p+1)}(t)=\sum_{k\in\mathbb{Z}^n}\widehat{\mathbf{U}}^{(p+1)}_k(t)e^{ik\cdot x}.
\end{align*}
Applying the mean value theorem, we get
\begin{align*}
    T_h:&=\left\|\frac{\mathbf{U}^{(p)}(t+h)-\mathbf{U}^{(p)}(t)}{h}-\sum_{k\in\mathbb{Z}^n}\widehat{\mathbf{U}}^{(p+1)}_k(t)e^{ik\cdot x}\right\|^2_{\mathcal{H}^q(\T^n)}\\
    &=\sum_{k\in\mathbb{Z}^n}(1+\|k\|^2)^q\left|\frac{\widehat{\mathbf{U}}^{(p)}_k(t+h)-\widehat{\mathbf{U}}^{(p)}_k(t)}{h}-\widehat{\mathbf{U}}^{p+1}_k(t)\right|^2\\
    &=\left|\frac{\widehat{\mathbf{U}}^{(p)}_0(t+h)-\widehat{\mathbf{U}}^{(p)}_0(t)}{h}-\widehat{\mathbf{U}}^{p+1}_0(t)\right|^2\\
    &+\sum_{0\ne k\in\mathbb{Z}^n}(1+\|k\|^2)^q\left\|\cos\left(\sqrt{-M_k}(t+\xi_{p,k})+\frac{(p+1)\pi}{2}\right)\right.\\
    &\qquad \qquad \left. -\cos\left(\sqrt{-M_k}t+\frac{(p+1)\pi}{2}\right)\right\|^2\|\sqrt{-M_k}\|^{2(p-1)}\|\widehat{\mathbf{b}}_k\|^2
\end{align*}
for some $\xi_{p,k}\in(0,h)$. This shows that as $h$ goes to $0$, each term in the summation approaches $0$. To pass the limit inside the summation, we only need to show the uniform boundedness of
\begin{align*}
    a_k:=(1+\|k\|^2)^{q-S}\|\sqrt{-M_k}\|^{2(p-1)}
\end{align*}
Three respective cases $\beta<n$, $\beta=n$ and $\beta>n$ are now taken into consideration.\\
When $\beta<n$, Theorem~\ref{thm:asymptotics_2} and Theorem~\ref{thm:asymptotics_1} show that $\sqrt{-M_k}$ is bounded for large $\|k\|$. This gives the uniform boundedness of $a_k$ is established as $q\le S$.\\
When $\beta=n$, fix an $\epsilon>0$. Theorem~\ref{thm:asymptotics_2} and Theorem~\ref{thm:asymptotics_1} implies that there exists $C,r>0$ such that $\|\sqrt{-M_k}\|^2\le C\log\|k\|\le C(1+\|k\|^2)^{\epsilon/2}$ whenever $\|k\|\ge r$. Then for large $\|k\|$,
\begin{align*}
    a_k\le C(1+\|k\|^2)^{q-S+(p-1)\epsilon/2}.
\end{align*}
This yields the uniform boundedness of $a_k$ as long as $q-S+(p-1)\epsilon/2\le 0$. So given $q<S$, we can always pick approriate $\epsilon$ to satisfy this.\\
When $\beta>n$, from Theorem~\ref{thm:asymptotics_2} and Theorem~\ref{thm:asymptotics_1}, there exists $C,r>0$ such that $\|\sqrt{-M_k}\|^2\le C\|k\|^{\beta-n}\le C(1+\|k\|^2)^{(\beta-n)/2}$ whenever $\|k\|\ge r$. Then 
\begin{align*}
    a_k\le C(1+\|k\|^2)^{q-S+(p-1)\frac{\beta-n}{2}}
\end{align*}
which is uniformly bounded provided that $q-S+(p-1)\frac{\beta-n}{2}\le 0$.
\end{proof}

The next theorem is a direct application of Theorem~\ref{thm:for_1} and Theorem~\ref{thm:for_regularity}. The assumptions are to guarantee that the second derivative exists.
\begin{theorem}
Let $n\ge1$, $\delta>0$ and $\beta<n+2$. Let $\mathbf{b}\in \mathcal{H}^{S}(\T^n)$ and suppose that
\begin{enumerate}
    \item $q\le S$ in the case $\beta<n$,
    \item $q<S$ in the case $\beta=n$,
    \item $q\le S$ in the case $\beta>n$.
\end{enumerate}
Then the map $\mathbf{U}(t)$ defined by \eqref{eq:forcing} is the unique solution in the space $\mathcal{H}^q(\T^n)$ of the peridynamic equation
\begin{align}\label{eq:Nonhomo-peri-new}
   \begin{cases}
       \frac{d^2}{dt^2}\mathbf{U}(t)&=\mathcal{L}^{\delta,\beta} \mathbf{U}(t)+\mathbf{b},\\
    \mathbf{U}(0)&=0,\\
    \frac{d}{dt}\mathbf{U}(0)&=0.
   \end{cases}
\end{align}
Moreover, when $\beta\le n$, then $\mathbf{U}(t)\in C^\infty([0,\infty); \mathcal{H}^{S}(\T^n))$  and when $\beta>n$, \\then
$\mathbf{U}(t)\in C^2([0,\infty); \mathcal{H}^{S+\beta-n}(\T^n))$.
\end{theorem}

\subsection{Convergence to the Navier-Cauchy equation}\label{sec:Conv-forcing}
In this section, we will provide some results on the convergence of the solutions of the peridynamic equation \eqref{eq:Nonhomo-peri-new} to the  solution of the corresponding Navier-Cauchy equation for two types of  limits as $\delta \to 0^+$, with $\beta< n+2$ is being fixed,  or as $\beta\to (n+2)^-$, with $\delta>0$ is being fixed. The Navier-Cauchy equation is given by
\begin{align}\label{eq:Nonhomo-Navier}
    \begin{cases}
        \frac{d^2}{dt^2}\mathbf{U}^\mathcal{N}(t)&=\mathcal{N}\mathbf{U}^{\mathcal{N}}(t)+\mathbf{b},\; t>0\\
    \mathbf{U}^{\mathcal{N}}(0)&=0,\\
    \frac{d}{dt}\mathbf{U}^{\mathcal{N}}(0)&=0,
    \end{cases}
\end{align}
and its solution is given by
\begin{align}
    \mathbf{U}^{\mathcal{N}}(t):=
    \widehat{\mathbf{b}}_0\frac{t^2}{2}+\sum_{0\ne k\in \mathbb{Z}^n}\left[\cos\left(\sqrt{-M^{\mathcal{N}}_k}t\right)-I\right](M^{\mathcal{N}}_k)^{-1}\widehat{\mathbf{b}}_k e^{ik\cdot x}.
\end{align}
To emphasize the dependence of the solution of the nonlocal wave equation \eqref{eq:Nonhomo-peri-new} on the parameters $\delta$ and $\beta$, in the remaining part of this section we denote the solution $\mathbf{U}$ by $\mathbf{U}^{\delta,\beta}$.

\begin{theorem}\label{thm:forcing}
Let $n\ge 1$, $\delta>0$ and $\beta<n+2$. Suppose $\mathbf{b}\in \mathcal{H}^{S}(\T^n)$. Let $\mathbf{U}^{\delta,\beta}(t)\in \mathcal{H}^{s'}(\T^n)$ and  $\mathbf{U}^{\mathcal{N}}(t)\in \mathcal{H}^s(\T^n)$ be the solutions to the peridynamic equation \eqref{eq:Nonhomo-peri-new} and the Navier-Cauchy equation \eqref{eq:Nonhomo-Navier}, respectively, where $s'=S+\max\{0,\beta-n\}$ and $s=S+2$ . Let $t\ge 0$, when $\beta<n+2$ be fixed, then
\begin{align*}
    \lim_{\delta\to 0^+}\mathbf{U}^{\delta,\beta}(t)=\mathbf{U}^{\mathcal{N}}(t),\quad \text{in }\mathcal{H}^{s'}(\T^n).
\end{align*}
\end{theorem}
\begin{proof}
From \ref{eq:forcing}, one can check that
\begin{align*}
    \|\mathbf{U}(t)-\mathbf{U}^{\mathcal{N}}(t)\|^2_{\mathcal{H}^{s'}(\T^n)}\le\sum_{0\ne k\in \mathbb{Z}^n}(1+\|k\|^2)^{s'}&\left\|[\cos(\sqrt{-M_k}t)-I](M_k)^{-1}\right.\\
    &\qquad \left.-[\cos(\sqrt{-M^{\mathcal{N}}_k}t)-I](M^{\mathcal{N}}_k)^{-1}\right\|^2\|\widehat{\mathbf{b}}_k\|^2.
\end{align*}
In this summation, each term converges to $0$ using Proposition \ref{prop1}. Since $\mathbf{b}\in H^S(\T^n)$, we need the uniform boundedness of
\begin{align*}
    \|k\|^{s'-S}\left\|[\cos(\sqrt{-M_k}t)-I](M_k)^{-1}-[\cos(\sqrt{-M^{\mathcal{N}}_k}t)-I](M^{\mathcal{N}}_k)^{-1}\right\|
\end{align*}
for large $\|k\|$ to pass the limit inside the summation. Observe that
\begin{align*}
    &\left\|[\cos(\sqrt{-M_k}t)-I](M_k)^{-1}-[\cos(\sqrt{-M^{\mathcal{N}}_k}t)-I](M^{\mathcal{N}}_k)^{-1}\right\|\\
    &=\max\limits_{i=1,2}\left|\frac{\cos(\sqrt{\lambda_i(k)})-1}{-\lambda_i(k)}-\frac{\cos(\sqrt{\lambda^{\mathcal{N}}_i(k)})-1}{-\lambda^{\mathcal{N}}_i(k)}\right|\\
    &\le\max\limits_{i=1,2}\left(\frac{2}{|\lambda_i(k)|}+\frac{2}{|\lambda^{\mathcal{N}}_i(k)|}\right).
\end{align*}
From Theorem \ref{thm:asymptotics_2} and Theorem \ref{thm:asymptotics_1}, the following bound holds
\begin{align*}
    |\lambda_1(k)|,|\lambda_2(k)|\ge\begin{cases}
    C_1&\text{if }\beta\le n\\
    C_2\|k\|^{\beta-n}&\text{if }\beta>n
    \end{cases},
\end{align*}
whenver $\|k\|\ge r$ for positive constants $r, C_1, C_2$. This shows that
\begin{align*}
    \left\|[\cos(\sqrt{-M_k}t)-I](M_k)^{-1}-[\cos(\sqrt{-M^{\mathcal{N}}_k}t)-I](M^{\mathcal{N}}_k)^{-1}\right\|\le\begin{cases}
    \frac{D}{\|k\|^2}+C&\text{if }\beta\le n\\
    \frac{D}{\|k\|^2}+\frac{C}{\|k\|^{\beta-n}}&\text{if }\beta> n
    \end{cases}
\end{align*}
Moreover,
\begin{align*}
    \|k\|^{s'-S}\le\begin{cases}
    1&\text{if }\beta\le n\\
    \|k\|^{\beta-n}&\text{if }\beta>n
    \end{cases}
\end{align*}
These properties give the required uniform boundedness.
\end{proof}

\begin{theorem}
Let $n\ge 1$, $\delta>0$ and $\beta<n+2$. Suppose $\mathbf{b}\in \mathcal{H}^{S}(\T^n)$. Let $\mathbf{U}^{\delta,\beta}(t)\in \mathcal{H}^{s'}(\T^n)$ and  $\mathbf{U}^{\mathcal{N}}(t)\in \mathcal{H}^s(\T^n)$ be the solutions to the peridynamic equation \eqref{eq:Nonhomo-peri-new} and the Navier-Cauchy equation \eqref{eq:Nonhomo-Navier}, respectively, where $s'=S+\max\{0,\beta-n\}$ and $s=S+2$ . Let $t\ge 0$, then, for $0<\epsilon<2$,
\begin{align*}
    \lim_{\beta\to n+2^-}\mathbf{U}^{\delta,\beta}(t)=\mathbf{U}^{\mathcal{N}}(t),\quad \text{in }\mathcal{H}^{s_0}(\T^n),
\end{align*}
where $s_0=S+2-\epsilon$.
\end{theorem}
\begin{proof}
With fixed $\epsilon\in (0,2)$, for any $\beta<n+2$ such that $n+2-\beta<\epsilon$, both $\mathcal{H}^s(\T^n), \mathcal{H}^{s'}(\T^n)$ are subspaces of $\mathcal{H}^{s_0}(\T^n)$, hence the limit makes sense. Notice that
\begin{align*}
    \|\mathbf{U}(t)-\mathbf{U}^{\mathcal{N}}(t)\|^2_{\mathcal{H}^{s'}(\T^n)}\le\sum_{0\ne k\in \mathbb{Z}^n}(1+\|k\|^2)^{s'}&\left\|[\cos(\sqrt{-M_k}t)-I](M_k)^{-1}\right.\\
    &\qquad\left.    -[\cos(\sqrt{-M^{\mathcal{N}}_k}t)-I](M^{\mathcal{N}}_k)^{-1}\right\|^2\|\widehat{\mathbf{b}}_k\|^2.
\end{align*}
The proof is complete, provided the uniform boundedness of
\begin{align*}
    \|k\|^{2-\epsilon}\left\|[\cos(\sqrt{-M_k}t)-I](M_k)^{-1}-[\cos(\sqrt{-M^{\mathcal{N}}_k}t)-I](M^{\mathcal{N}}_k)^{-1}\right\|
\end{align*}
as $\|k\|\to \infty$. Since $\beta>n$, in Theorem \ref{thm:forcing}, we have proved that
\begin{align*}
    \|k\|^{2-\epsilon}\left\|[\cos(\sqrt{-M_k}t)-I](M_k)^{-1}-[\cos(\sqrt{-M^{\mathcal{N}}_k}t)-I](M^{\mathcal{N}}_k)^{-1}\right\|&\le \|k\|^{2-\epsilon}\left(\frac{D}{\|k\|^2}+\frac{C}{\|k\|^{\beta-n}}\right)\\
    &=\frac{D}{\|k\|^\epsilon}+\frac{C}{\|k\|^{\epsilon+\beta-(n+2)}},
\end{align*}
which is obviously bounded.
\end{proof}

\subsection{Spatial and temporal regularity over \texorpdfstring{$(L^2(\T^n))^n$}{}}\label{sec:temp-forcing}
The solution $\mathbf{U}(t)\in \mathcal{H}^s(\T^n)$ given in \eqref{eq:forcing} is a distribution when $s<0$, and  when $s\ge 0$ defines a regular vector field
\begin{align*}
     \mathbf{u}(x,t)=\mathbf{U}(t)(x)
    =\widehat{\mathbf{b}}_0\frac{t^2}{2}+\sum_{0\ne k\in \mathbb{Z}^n}\left[\cos\left(\sqrt{-M_k}t\right)-I\right]M_k^{-1}\widehat{\mathbf{b}}_k e^{ik\cdot x}.
\end{align*}
Let the data be a regular vector field; $\mathbf{b}\in \mathcal{H}^{S}(\T^n)$ with $S\ge 0$. In this section, we focus on conditions that guarantee that $\mathbf{u}(x,t)$, the solution of \eqref{eq:Nonhomo-peri}, is a regular vector field $\mathbf{u}(\cdot,t)\in (L^2(\T^n))^n$.


\subsubsection{Temporal regularity with respect to the Gateaux derivative}
\begin{theorem}
Let $n\ge 1$, $\delta>0$ and $\beta<n+2$. Suppose that
\begin{enumerate}
    \item $s\ge 0$ in the case when $\beta<n$,
    \item $s>0$ in the case when $\beta=n$, and
    \item $s\ge 0 $ in the case $\beta>n$.
\end{enumerate}
Then the  vector field $\mathbf{u}$, where $  \mathbf{u}(x,t)=\mathbf{U}(t)(x)$, with $\mathbf{U}$ defined  in \eqref{eq:forcing},  is the unique solution to the  peridynamic equation \eqref{eq:Nonhomo-peri}.
Moreover, when $\beta\le n$, then $\mathbf{u}(x,t)\in C^\infty([0,\infty); \mathcal{H}^S(\T^n))$, and when $\beta>n$, then
$\mathbf{u}(x,t)\in C^2([0,\infty); \mathcal{H}^{S+\beta-n}(\T^n))$.
\end{theorem}

\subsubsection{Temporal regularity with respect to the classical derivative}
\begin{theorem}
    Let $n\ge 1$, $\delta>0$ and $\beta<n$. Suppose that $\mathbf{b}$ satisfies the $\alpha$ - H{\"o}lder continuity condition with $\alpha>n/2$. Then $\mathbf{u}(x,\cdot)\in C^\infty[0,\infty)$, for any $x\in\T^n$.
\end{theorem}

\begin{proof}
    Fix a value of $p$ in $\mathbb{Z}_{\ge 0}$. Consider the series 
    $$\displaystyle \sum_{0\ne k\in \mathbb{Z}^n}\left[\cos\left(\sqrt{-M_k}t\right)-I\right]M_k^{-1}\widehat{\mathbf{b}}_k e^{ik\cdot x}.$$ From Theorem \ref{thm:asymptotics_2} and Theorem \ref{thm:asymptotics_1}, $\|\sqrt{-M_k}\|$ and $\|({-M_k})^{-1}\|$ are bounded for large $k$. This makes
\begin{align*}
    &\left\|\left[\cos\left(\sqrt{-M_k}t\right)-I\right]M_k^{-1}\widehat{\mathbf{b}}_k e^{ik\cdot x}\right\|\le2\|M_k^{-1}\|\|\widehat{\mathbf{b}}_k\|\le C_1\|\widehat{\mathbf{b}}_k\|\\
    &\left\|\frac{d^p}{dt^p}\left[\cos\left(\sqrt{-M_k}t\right)-I\right]M_k^{-1}\widehat{\mathbf{b}}_k e^{ik\cdot x}\right\|\le\|\sqrt{-M_k}\|^{p-2}\|\widehat{\mathbf{b}}_k\| \le C_2\|\widehat{\mathbf{b}}_k\|
\end{align*}
as $k\to \infty$. From Lemma \ref{Fourier_summable}, the Fourier coefficients of $\widehat{\mathbf{b}}_k$ is absolutely summable. The proof is  complete by using Lemma \ref {lem:series}.
\end{proof}

\begin{theorem}
    Let $n\ge 1$, $\delta>0$ and $\beta=n$. Suppose $S>n$. Then $\mathbf{u}(x,\cdot)\in C^\infty[0,\infty)$, for any $x\in\T^n$. 
\end{theorem}

\begin{proof}
    Fix a value of $p$ in $\mathbb{Z}_{\ge 0}$ and choose $\epsilon>0$ such that $S-(p-2)\epsilon>n$. Consider the series $$\displaystyle \sum_{0\ne k\in \mathbb{Z}^n}\left[\cos\left(\sqrt{-M_k}t\right)-I\right]M_k^{-1}\widehat{\mathbf{b}}_k e^{ik\cdot x}.$$ From Theorem \ref{thm:asymptotics_2} and Theorem \ref{thm:asymptotics_1}, 
    \begin{align*}
    \|\sqrt{-M_k}\|\le C(\log\|k\|)^{1/2}\le C\|k\|^{\epsilon},\quad\text{and}\quad
    \|({-M_k})^{-1}\|\le  C,
\end{align*}
for large $k$. This implies
\begin{align*}
    &\left\|\left[\cos\left(\sqrt{-M_k}t\right)-I\right]M_k^{-1}\widehat{\mathbf{b}}_k e^{ik\cdot x}\right\|\le2\|M_k^{-1}\|\|\widehat{\mathbf{b}}_k\|\le \frac{C_1}{\|k\|^S}\\
    &\left\|\frac{d^p}{dt^p}\left[\cos\left(\sqrt{-M_k}t\right)-I\right]M_k^{-1}\widehat{\mathbf{b}}_k e^{ik\cdot x}\right\|\le\|\sqrt{-M_k}\|^{p-2}\|\widehat{\mathbf{b}}_k\| \le \frac{C_2}{\|k\|^{S-(p-2)\epsilon}}
\end{align*}
as $k\to \infty$. The proof is  complete by using Lemma \ref {lem:series}.
\end{proof}

\begin{theorem}
    Let $n\ge 1$, $\delta>0$ and $n<\beta<n+2$. Suppose that there exist a nonnegative integer $p$ such that 
    \begin{align*}
        S>n+(p-2)\frac{\beta-n}{2}.
    \end{align*} Then $\mathbf{u}(x,\cdot)\in C^p[0,\infty)$, for any $x\in\T^n$.
\end{theorem}

\begin{proof}
    Fix a value of $p$ in $\mathbb{Z}_{\ge 0}$. Consider the series $$\displaystyle \sum_{0\ne k\in \mathbb{Z}^n}\left[\cos\left(\sqrt{-M_k}t\right)-I\right]M_k^{-1}\widehat{\mathbf{b}}_k e^{ik\cdot x}.$$ From Theorem \ref{thm:asymptotics_2} and Theorem \ref{thm:asymptotics_1}, 
    \begin{align*}
    \|\sqrt{-M_k}\|\le C\|k\|^{(\beta-n)/2},\quad\text{and}\quad
    \|({-M_k})^{-1}\|\le  C\|k\|^{-(\beta-n)},
\end{align*}
for large $k$. This implies
\begin{align*}
    &\left\|\left[\cos\left(\sqrt{-M_k}t\right)-I\right]M_k^{-1}\widehat{\mathbf{b}}_k e^{ik\cdot x}\right\|\le2\|M_k^{-1}\|\|\widehat{\mathbf{b}}_k\|\le \frac{C_1}{\|k\|^{S+\beta-n}}\\
    &\left\|\frac{d^p}{dt^p}\left[\cos\left(\sqrt{-M_k}t\right)-I\right]M_k^{-1}\widehat{\mathbf{b}}_k e^{ik\cdot x}\right\|\le\|\sqrt{-M_k}\|^{p-2}\|\widehat{\mathbf{b}}_k\| \le \frac{C_2}{\|k\|^{S-(p-2)(\beta-n)/2}}
\end{align*}
as $k\to \infty$. The proof is  complete using Lemma \ref {lem:series}.
\end{proof}

\clearpage
\bibliographystyle{acm}
\bibliography{ref}

\begin{thebibliography}{10}

\bibitem{alali2021fourier}
{\sc Alali, B., and Albin, N.}
\newblock Fourier multipliers for nonlocal laplace operators.
\newblock {\em Applicable Analysis 100}, 12 (2021), 2526--2546.

\bibitem{alali2022}
{\sc Alali, B., and Albin, N.}
\newblock Linear peridynamics fourier multipliers and eigenvalues.
\newblock {\em Journal of Peridynamics and Nonlocal Modeling 6}, 2 (2024),
  294--317.

\bibitem{alali2012multiscale}
{\sc Alali, B., and Lipton, R.}
\newblock Multiscale dynamics of heterogeneous media in the peridynamic
  formulation.
\newblock {\em Journal of Elasticity 106\/} (2012), 71--103.

\bibitem{alimov2019solvability}
{\sc Alimov, S., and Sheraliev, S.}
\newblock On the solvability of the singular equation of peridynamics.
\newblock {\em Complex Variables and Elliptic Equations 64}, 5 (2019),
  873--887.

\bibitem{alimov2020solvability}
{\sc Alimov, S., and Sheraliev, S.}
\newblock On the solvability of hypersingular equation of peridynamics.
\newblock {\em Bulletin of National University of Uzbekistan: Mathematics and
  Natural Sciences 3}, 3 (2020), 278--298.

\bibitem{alimov2014problems}
{\sc Alimov, S.~A., Cao, Y., and Ilhan, O.}
\newblock On the problems of peridynamics with special convolution kernels.

\bibitem{alimov2023hypersingular}
{\sc Alimov, S.~A., and Sheraliev, S.~N.}
\newblock On hypersingular operators associated with peridynamics.
\newblock {\em Differential Equations 59}, 7 (2023), 914--918.

\bibitem{dang2024}
{\sc Dang, T., Alali, B., and Albin, N.}
\newblock Regularity of solutions for the nonlocal wave equation on periodic
  distributions.
\newblock {\em Submitted\/} (2024).

\bibitem{NIST:DLMF}
{\it NIST Digital Library of Mathematical Functions}, Release 1.0.17 of
  2017-12-22.
\newblock F.~W.~J. Olver, A.~B. {Olde Daalhuis}, D.~W. Lozier, B.~I. Schneider,
  R.~F. Boisvert, C.~W. Clark, B.~R. Miller and B.~V. Saunders, eds.

\bibitem{emmrichperidynamic}
{\sc Emmrich, E., and Weckner, O.}
\newblock The peridynamic equation of motion in non-local elasticity theory.
\newblock In {\em III European Conference on Computational Mechanics: Solids,
  Structures and Coupled Problems in Engineering: Book of Abstracts}, Springer,
  pp.~62--62.

\bibitem{Folland1999}
{\sc Folland, G.}
\newblock {\em Real analysis: modern techniques and their applications},
  vol.~40.
\newblock John Wiley \& Sons, 1999.

\bibitem{foss2016differentiability}
{\sc Foss, M., and Radu, P.}
\newblock Differentiability and integrability properties for solutions to
  nonlocal equations.
\newblock In {\em New Trends in Differential Equations, Control Theory and
  Optimization: Proceedings of the 8th Congress of Romanian Mathematicians\/}
  (2016), World Scientific, pp.~105--119.

\bibitem{foss2018existence}
{\sc Foss, M.~D., Radu, P., and Wright, C.}
\newblock Existence and regularity of minimizers for nonlocal energy
  functionals.
\newblock {\em Differential and Integral Equations 31}, 11-12 (2018), 807--832.

\bibitem{Grafakos2009}
{\sc Grafakos, L.}
\newblock {\em Classiccal Fourier Analysis}.
\newblock Springer, 2008.

\bibitem{mengesha2014nonlocal}
{\sc Mengesha, T., and Du, Q.}
\newblock Nonlocal constrained value problems for a linear peridynamic navier
  equation.
\newblock {\em Journal of Elasticity 116}, 1 (2014), 27--51.

\bibitem{ilyas}
{\sc Mustapha, I., Alali, B., and Albin, N.}
\newblock Regularity of solutions for nonlocal diffusion equations on periodic
  distributions.
\newblock {\em Journal of Integral Equations and Applications 35}, 1 (2023),
  81--104.

\end{thebibliography}

\end{document}